\documentclass{article}
\usepackage[journal=APDE,lang=british]{ems-journal-jems} 

\usepackage{textcomp}
\usepackage{mathrsfs}

\usepackage{chngcntr}
\usepackage{apptools}
\usepackage{verbatim}
\usepackage{tikz-cd}
\usepackage{csquotes}

\newtheorem{theorem}{Theorem}
\newtheorem*{theorem*}{Theorem}

\newtheorem{proposition}{Proposition}
\newtheorem{lemma}{Lemma}

\theoremstyle{definition}

\theoremstyle{remark}
\newtheorem{remark}{Remark}
\newtheorem{example}{Example}
\newtheorem{corollary}{Corollary}

\numberwithin{equation}{section}

\newcommand\restr[2]{{
		\left.\kern-\nulldelimiterspace 
		#1 
		\vphantom{\big|} 
		\right|_{#2} 
}}

\numberwithin{equation}{section}

\numberwithin{theorem}{section}
\numberwithin{proposition}{section}
\numberwithin{lemma}{section}
\numberwithin{remark}{section}
\numberwithin{example}{section}
\numberwithin{exercise}{section}
\numberwithin{corollary}{section}

\newcounter{desccount}

\begin{document}

\title{On the Bendixson criterion for attractors of dynamical systems in metric spaces}
\titlemark{On the Bendixson criterion for attractors of dynamical systems in metric spaces}


\emsauthor{1}{
	\givenname{Mikhail}
	\surname{Anikushin}
	\orcid{0000-0002-7781-8551}}{M.M.~Anikushin}

\Emsaffil{1}{
	\pretext{}
	\department{Department of Applied Cybernetics}
	\organisation{Faculty of Mathematics and Mechanics, St Petersburg University}
	\rorid{023znxa73}
	\address{Universitetskiy prospekt 28}
	\zip{198504}
	\city{Peterhof}
	\country{Russia}
	\posttext{}
	\affemail{demolishka@gmail.com}
}
	

	


\classification[37C25, 37C60, 37L15, 37L30, 55N05]{37C45}

\begin{abstract}
We revisit the Bendixson criterion for attractors, originally proposed by R.A.~Smith for the finite-dimensional case and subsequently extended by M.Y.~Li and J.S.~Muldowney to Banach spaces. First, we give a concise proof of the criterion for general dynamical systems (semiflows, mappings, cocycles) in metric spaces using \v{C}ech cohomology. This approach allows for a significant weakening of previously known conditions while simultaneously resolving existing misunderstandings. Second, for compact invariant sets of semiflows that are topologically one-dimensional and possess zero first cohomology (e.g., attracting a superset with zero first cohomology), we prove that any chain recurrent point within the invariant set must be an equilibrium. In the proof, we construct a global strict Lyapunov function using Whitney--Bebutov sections. Third, we generalize and refine R.A.~Smith's theorem concerning the finiteness of the set of fixed points for analytic mappings, as well as M.L.~Cartwright's theorems on to almost periodic flows.
\end{abstract}

\maketitle

\section{Introduction}

For vector fields on the plane, the classical Bendixson negative criterion rules out the existence of closed invariant contours $\Gamma$ lying entirely within a simply connected domain $\mathcal{D}$, provided that the divergence maintains a constant sign in $\mathcal{D}$ \cite{Bendixson1901}. It appears that Smith was the first to note that this criterion naturally generalizes to higher-dimensional systems when viewed in terms of attractors \cite{Smith1986HD}. In this case, the relevant conditions are formulated in terms of dimension estimates \cite{KuzReit2020, ZelikAttractors2022, Anikushin2023LyapExp}. Specifically, Smith's approach requires the Hausdorff dimension $\dim_{\operatorname{H}}\mathcal{A}$ of an attractor $\mathcal{A}$ to be strictly less than $2$. Later, Li and Muldowney extended this criterion to maps and semiflows in Banach spaces to rule out the existence of certain $m$-dimensional bounding surfaces, given the condition $\dim_{\operatorname{H}}\mathcal{A} < m$ \cite{LiMuldowney1995LowBounds}. We discuss some technical nuances of these works at the end of Section \ref{SEC: CohomologicalBendixsonCriterion}.

Our primary interest in the topic is concerned with applications to global stability problems for semiflows. This also originates from the paper \cite{Smith1986HD}, where Smith noted that, in applications to ODEs, conditions of uniform area contraction --- which particularly guarantee the estimate $\dim_{\operatorname{H}}\mathcal{A} < 2$ --- are preserved under small $C^{1}$-perturbations of the system. Building on this, he employed Pugh's closing lemma to prove the global stability of the system in the strong sense\footnote{We follow the terminology of \cite{Kuzetal2020Lorenz}, according to which a system is called \textit{globally stable} if any of its trajectories tends to the stationary set (possibly infinite). If any trajectory tends to a single equilibrium, we call the semiflow \textit{globally stable in the strong sense}.} --- that is, the convergence of any trajectory to an equilibrium. Subsequently, this approach was popularized in \cite{LiMuldowney1996SIAMGlobStab}, and the conditions for uniform area contraction were refined using various metrics \cite{LeonovBoi1992, NoackReitmann1996, KuzReit2020, Kuzetal2020Lorenz}. A more detailed discussion on this can be found in our papers \cite{Anikushin2023LyapExp, AnikushinRomanov2023FreqConds}.

Recently, we obtained effective dimension estimates for delay equations using either explicitly constructed metrics \cite{Anikushin2023LyapExp, AnikushinRomanov2024EffEst} or implicit quadratic functionals on exterior powers \cite{AnikushinRomanovQBC2026,AnikushinRomanov2026Subroutine,AnikushinRomanov2023FreqConds}. In the latter case, the results are based on the development of the functional-analytical apparatus of frequency-domain methods presented in our works \cite{Anikushin2023Comp, Anikushin2020FreqDelay, Anikushin2020FreqParab}. These general methods are potentially applicable to a wide range of problems that exhibit some kind of asymptotic compactness. Although the frequency conditions arising from this approach can only be verified numerically, experiments in \cite{AnikushinRomanovQBC2026,AnikushinRomanov2026Subroutine,AnikushinRomanov2023FreqConds} demonstrate significant improvements over certain known results regarding global stability obtained using more specialized methods.

To the best of our knowledge, there are currently no versions of the closing lemma applicable to infinite-dimensional systems. Consequently, it is impossible to derive global stability based on area contraction using the classical line of reasoning. This prompted us to seek purely topological mechanisms ensuring global stability. In the present work, we explicate such mechanisms based on \v{C}ech cohomology and topological dynamics.

There are related developments known as the shape theory of attractors. These are based on the shape equivalence between attractors and their domains of attraction, enabling a topological analysis of bifurcations \cite{KapitanskiRodnianski2000, BargeGiraldoSanjurjo2020}. However, we will need only surjectivity in cohomology, which holds under much weaker assumptions. 

More precisely, in Section \ref{SEC: CohomologicalBendixsonCriterion}, we present a cohomological perspective on the Bendixson criterion for attractors, based on just two principles. First, the cohomology of an attractor $\mathcal{A}$ cannot be larger than cohomology of any superset that it attracts. Second, the constraint on the topological dimension $\dim_{\operatorname{T}}\mathcal{A} \leq m$ allows to carry this conclusion to the $m$th cohomology group of any compact subset $\Gamma$ of $\mathcal{A}$. This applies to general dynamical systems (including nonautonomous ones) with discrete or continuous time, although in the latter case the conditions are significantly less stringent and do not require injectivity, owing to the existence of a natural homotopy. Specific variants of the Bendixson criterion are stated in Theorems \ref{TH: BendixsonSemiflow}, \ref{TH: BendixsonCocycleContTime}, \ref{TH: BendixsonMapping}, and \ref{TH: BendixsonCocyclesDiscreteTime}.

In the special case $m=1$ for semiflows, we obtain the following more stronger result.
\begin{theorem}
	\label{TH: BendixsonGlobalStabilityOfSemiflows}
	Let $\varphi$ be a continuous semiflow in a metric space $\mathcal{X}$, and let $\mathcal{A}$ be a compact invariant subset: $\varphi^{t}(\mathcal{A}) = \mathcal{A}$ for any $t \geq 0$. Suppose that
	\begin{enumerate}
		\item[1)] $\check{H}^{1}(\mathcal{A};\mathbb{Z}) = 0$ (e.g., if $\mathcal{A}$ attracts a superset $\mathcal{B} \supset \mathcal{A}$ with $\check{H}^{1}(\mathcal{B};\mathbb{Z}) = 0$).
		\item[2)] $\dim_{\operatorname{T}}\mathcal{A} \leq 1$;
	\end{enumerate}
	Then any chain recurrent point of $(\mathcal{A}, \varphi)$ is an equilibrium. In particular, the following holds:
	\begin{enumerate}
		\item[(a)] For any nonstationary complete trajectory $\gamma$ of $(\mathcal{A}, \varphi)$, the $\alpha$-limit set $\alpha(\gamma)$ and the $\omega$-limit set $\omega(\gamma)$ of $\gamma$ are disjoint connected subsets of equilibria;
		\item[(b)] More generally, for any finite family of nonstationary complete trajectories $\gamma_{1},\ldots,\gamma_{n}$ in $\mathcal{A}$ such that $\omega(\gamma_{i}) \cap \alpha(\gamma_{i+1}) \not= \emptyset$ for any $i \in \{1,\ldots,n-1\}$, we necessarily have that $\omega(\gamma_{n}) \cap \alpha(\gamma_{1}) = \emptyset$.
	\end{enumerate}
\end{theorem}
When $\mathcal{A}$ attracts compact subsets of $\mathcal{X}$, any nonwandering point of $(\mathcal{X},\varphi)$ is a chain recurrent point of $(\mathcal{A},\varphi)$, see Lemma \ref{LEM: NonwanderingAreChainRecurrent}. If, in addition, $\mathcal{X}$ is a normed vector space, then $\mathcal{A}$ attracts its closed convex hull, which is contractible; thus, the condition from item 2) is satisfied. In this case, the conclusion of Theorem \ref{TH: BendixsonGlobalStabilityOfSemiflows} is even stronger than what is obtained from the closing lemma and robustness of dimension estimates in \cite{Smith1986HD, LiMuldowney1996SIAMGlobStab}. Furthermore, if the set of equilibria is finite, then we immediately obtain that the semiflow is globally stable in the strong sense. This is because the $\omega$-limit set of a point is connected and consists of nonwandering points (which are also chain recurrent within the $\omega$-limit set). By virtue of item (b), the attractor $\mathcal{A}$ consists of equilibria and complete orbits connecting them which cannot form a cycle. Since $\mathcal{A}$ is connected, the number of connecting orbits equals the number of equilibria minus one.

Without finiteness of equilibria, the semiflow is globally stable, though not necessarily in the strong sense. Indeed, it is not hard to construct a flow on the plane whose attractor is a segment of equilibria onto which trajectories wind in a spiral. In this example, the equilibria are degenerate in all directions. In practice, the condition $\dim_{\operatorname{T}}\mathcal{A} \leq 1$ is verified by means of the uniform contraction of infinitesimal areas over $\mathcal{A}$, see \eqref{EQ: BanachSpaceAreaContraction}. In particular, this guarantees that any equilibrium can have at most one degenerate direction. In this case, by invoking center manifold theory, one can establish global stability in the strong sense, even when the set of equilibrium states is infinite, see \cite[Theorem IV.1]{AnikushinRomanovQBC2026}.

A proof of Theorem \ref{TH: BendixsonGlobalStabilityOfSemiflows} is given in Section \ref{SEC: GlobalStabilitySemiflowsBendixson}. First, we reduce the problem to the case of flows by transiting to the natural extension $(\hat{\mathcal{A}},\hat{\varphi})$ of $(\mathcal{A},\varphi)$. This reduction is possible because cohomology, upper bounds on dimension, and chain recurrence behave well under taking inverse limits. Next, for $(\hat{\mathcal{A}},\hat{\varphi})$, we construct a global (strict) Lyapunov function glued from local Lyapunov functions associated with the Whitney--Bebutov sections, which are zero-dimensional in our case and thereby ensure continuity of the local functions. In our construction, the global Lyapunov function takes the set of equilibria to a zero measure subset, and this suffices to obtain the desired conclusion.

In Section \ref{SEC: IsolatedFixedPointsAnMaps}, we present applications of the Bendixson criterion for deriving conditions under which the set of fixed points of analytic mappings is finite. This generalizes the corresponding results of \cite{Smith1986HD, NakajimaSeifert1983} to infinite dimensions.

In Section \ref{SEC: AlmostPeriodicBendixsonRefinement}, to illustrate applications of the Bendixson criterion to nonautonomous systems, we generalize and refine certain results by Cartwright concerning almost periodic flows \cite{Cartwright1967, Cartwright1969}, making use of the theory of compact abelian groups. We then discuss the difficulties that arise when attempting to extend these results to the case of almost automorphic systems.
\section{Preliminaries}
\label{SEC: Preliminaries}
In this section, we introduce several concepts related to dynamical systems, primarily to establish notation and terminology. More specific concepts will be presented in the relevant sections.

Let $\mathcal{Q}$ be a metric space, and let $\mathbb{T} \in \{ \mathbb{Z}_{+}, \mathbb{Z}, \mathbb{R}_{+}, \mathbb{R} \}$ be a time set. A continuous dynamical system on $\mathcal{Q}$ is a family of mappings $\vartheta^{t} \colon \mathcal{Q} \to \mathcal{Q}$, where $t \in \mathbb{T}$, satisfying the following properties:
\begin{enumerate}
	\item[1)] $\vartheta^{0}(q) = q$ and $\vartheta^{t+s}(q) = \vartheta^{t}(\vartheta^{s}(q))$ for any $q \in \mathcal{Q}$ and $t,s \in \mathbb{T}$;
	\item[2)] the mapping $\mathbb{T} \times \mathcal{Q} \ni (t,q) \mapsto \vartheta^{t}(x) \in \mathcal{Q}$ is continuous.
\end{enumerate}
In this context, we write $(\mathcal{Q}, \mathbb{T}, \vartheta)$, $(\mathcal{Q},\vartheta)$, or simply $\vartheta$ to denote the dynamical system. Any $\vartheta^{t}$ is called the time-$t$ mapping of $\vartheta$.  When $T \in \{\mathbb{R}_{+}, \mathbb{R}\}$ (resp. $T \in \{\mathbb{Z}_{+}, \mathbb{Z}\}$), $\vartheta$ is said to be a system with \textit{continuous} (resp. \textit{discrete}) time. If $\mathbb{T} \in \{ \mathbb{Z}, \mathbb{R} \}$, then $\vartheta$ is said to be \textit{reversible} (or with \textit{reversible time}). Set also $\mathbb{T}^{+} \coloneq \mathbb{T} \cap [0,+\infty)$ and $\overline{\mathbb{T}} \coloneq \mathbb{T} \cup -\mathbb{T}$.

If $\mathbb{T} = \mathbb{R}_{+}$ (resp. $\mathbb{T}=\mathbb{R}$), then $\vartheta$ is called a continuous \textit{semiflow} (resp. \textit{flow}) on $\mathcal{X}$. For $\mathbb{T} \in \{\mathbb{Z}_{+}, \mathbb{Z}\}$, the entire family is determined by the time-$1$ mapping $\vartheta^{1}$, which will be denoted just by $T$.

Given $(\mathcal{Q},\vartheta)$, a subset $\mathcal{A} \subset \mathcal{Q}$ is said to \textit{attract} another subset $\mathcal{B} \subset \mathcal{Q}$ if for any open neighborhood $\mathcal{U}$ of $\mathcal{A}$ there exists $t_{0} \in \mathbb{T}^{+}$ such that $\vartheta^{t}(\mathcal{B}) \subset \mathcal{U}$ for any $t \geq t_{0}$.

Recall that a \textit{complete trajectory} of $\vartheta$ is a mapping $\gamma \colon \overline{\mathbb{T}} \to \mathcal{Q}$ such that $\vartheta^{t}(\gamma(s)) = \gamma(t+s)$ for any $t \in \mathbb{T}$ and $s \in \overline{\mathbb{T}}$.

Let $(\mathcal{X}, \mathbb{T}_{\pi}, \pi)$ and $(\mathcal{Q}, \mathbb{T}_{\vartheta}, \vartheta)$ be two dynamical systems, both with either discrete or continuous time. A surjective continuous mapping $\wp \colon \mathcal{X} \to \mathcal{Q}$ is called an \textit{extension} if
\begin{equation}
	\wp \circ \pi^{t} = \vartheta^{t} \circ \wp \qquad \text{for any} \ t \in \mathbb{T}_{\pi} \cap \mathbb{T}_{\vartheta}.
\end{equation}
For brevity, this situation is denoted as $\wp \colon (\mathcal{X}, \pi) \to (\mathcal{Q}, \vartheta)$. One also says that $(\mathcal{X}, \pi)$ is an extension of $(\mathcal{Q}, \vartheta)$ via $\wp$, and that $\wp$ (or the entire triple) is an extension of $(\mathcal{Q}, \vartheta)$ to $(\mathcal{X}, \pi)$. In this context, for any $q \in \mathcal{Q}$, the set $\mathcal{X}_{q} \coloneq \wp^{-1}(q)$ is called the \textit{fiber} over $q$. 

Moreover, when $\mathbb{T}_{\vartheta} \in \{ \mathbb{R}, \mathbb{Z} \}$ and $q \in \mathcal{Q}$ is fixed, a subset $\mathcal{A} \subset \mathcal{X}_{q}$ is said to \textit{pullback attract} over $q$ a family of subsets $\{\mathcal{B}_{s}\}_{s \in \mathbb{T}^{-}_{\vartheta}}$, where $\mathbb{T}^{-}_{\vartheta} \coloneq \mathbb{T}_{\vartheta} \cap (-\infty, 0]$ and $\mathcal{B}_{s} \subset \mathcal{X}_{\vartheta^{s}(q)}$ for any $s \in \mathbb{T}^{-}_{\vartheta}$, if for any open neighborhood $\mathcal{U}$ of $\mathcal{A}$ there exists $t_{0} \in \mathbb{T}^{+}_{\pi}$ such that $\pi^{t}(\mathcal{B}_{-t}) \subset \mathcal{U}$ for any $t \geq t_{0}$.

Let $\wp_{n} \colon (\mathcal{X}_{n+1}, \pi_{n+1}) \to (\mathcal{X}_{n}, \pi_{n})$, where $n=0,1,\ldots$, be a sequence of extensions of dynamical systems having the same time set $\mathbb{T}$. The inverse limit of this sequence,
\begin{equation}
	(\mathcal{X},\pi) = \varprojlim (\mathcal{X}_{n+1}, \pi_{n+1}),
\end{equation}
is a dynamical system $\pi$ defined on the set $\mathcal{X}$ consisting of all $x=(x_{0},x_{1},\ldots)$, where $x_{n} \in \mathcal{X}_{n}$ and $\wp_{n}(x_{n+1}) = x_{n}$ for any $n=0,1,\ldots$, by the rule
\begin{equation}
	\pi^{t}(x) = (\pi^{t}_{1}(x_{1}), \pi^{t}_{2}(x_{2}),\ldots) \qquad \text{for any} \ t \in \mathbb{T} \ \text{and} \ x \in \mathcal{X}.
\end{equation}
Since each $\mathcal{X}_{n}$ is a metric space, one can define a metric on $\mathcal{X}$ compatible with the subspace topology in the product $\prod_{n=0}^{\infty}\mathcal{X}_{n}$. Then $(\mathcal{X},\pi)$ is seen to be a continuous dynamical system that is an extension of any $(\mathcal{X}_{n},\pi_{n})$ via the projection $x=(x_{0},x_{1},\ldots) \mapsto x_{n}$.

Let $(\mathcal{Y},\varphi)$ be a dynamical system with time set $\mathbb{T} \in \{\mathbb{Z}_{+}, \mathbb{R}_{+}\}$ such that any time-$t$ mapping $\varphi^{t}$ is surjective. Then we can consider a sequence of extensions  $(\mathcal{X}_{n},\pi_{n})$ as above with $\mathcal{X}_{n} \coloneq \mathcal{X}$, $\pi_{n} = \varphi$, and $\wp_{n} \coloneq \varphi^{1}$ for any $n=0,1,\ldots$. The inverse limit $(\hat{\mathcal{X}},\hat{\varphi})$ of this sequence is called the \textit{natural extension} of $\varphi$. One can extend $\hat{\varphi}$ to a reversible system on $\hat{\mathcal{X}}$. Namely, let $\sigma \colon \hat{\mathcal{X}} \to \hat{\mathcal{X}}$ be the shift map: $\sigma(x) = (x_{1},x_{2},\ldots)$ for $x=(x_{0},x_{1},\ldots)$. Then for any $x \in \hat{\mathcal{X}}$ and $t = k + s$ with nonnegative integer $k$ and $s \in [0,1] \cap \mathbb{T}$, we define 
\begin{equation}
	\hat{\varphi}^{-t}(x) \coloneq \sigma^{k}(\hat{\varphi}^{1-s}(x)).
\end{equation}
In fact, it is easy to see that $(\hat{\mathcal{X}}, \hat{\varphi})$ is topologicaly conjugate to the translation flow on the space of complete trajectories of $(\mathcal{X},\varphi)$ endowed with the topology of pointwise convergence. Nevertheless, the category-theoretical approach is helpful to lift certain properties of $(\mathcal{X}_{n},\pi_{n})$ to the inverse limit from a general point of view.
\section{A cohomological perspective on the Bendixson criterion for attractors}
\label{SEC: CohomologicalBendixsonCriterion}
During this section, $\mathcal{X}$ denotes a fixed metric space, and $G$ denotes an abelian group. For a metric space $\mathcal{Y}$, we will need the \v{C}ech cohomology modules $\check{H}^{*}(\mathcal{Y};G)$ of $\mathcal{Y}$ with coefficients in $G$ and certain related properties, which we are going to formulate. We assume that the reader is familiar with elements of algebraic topology and the related category-theoretical terminology, see \cite{Spanier1966, EilenbergSteenrod2015}. In particular, we assume familiarity with the fact that for any continuous mapping $f \colon \mathcal{X} \to \mathcal{Y}$, there is the associated homomorphism $f^{*} \colon \check{H}^{*}(\mathcal{X};G) \to \check{H}^{*}(\mathcal{Y};G)$ of cohomology modules, and that the correspondence $\mathcal{Y} \mapsto \check{H}^{*}(\mathcal{Y};G)$ is a contravariant functor --- that is, $(f \circ g)^{*} = g^{*} \circ f^{*}$ is satisfied for appropriate mappings, and if $f$ is the identity, then so is $f^{*}$.

Since every metric space is paracompact, the \v{C}ech cohomology is functorially isomorphic to the Alexander--Spanier cohomology, see \cite[Corollary 8, Section 6.8]{Spanier1966}\footnote{We emphasize that the referenced theorem establishes the isomorphism for cohomologies of a single space. For pairs of spaces, the related statement is left as an exercise, see \cite[Chapter 6, Exercises D1--D3]{Spanier1966}.}. In particular, \cite[Theorem 6, Section 6.5]{Spanier1966} gives the following:
\begin{description}[before=\let\makelabel\descriptionlabel]
	\item[\textbf{(CH1)}\refstepcounter{desccount}\label{DESC: CH1}] Homotopic mappings induce the same morphisms in cohomology modules.
\end{description}

Moreover, \cite[Theorem 2, Section 6.6]{Spanier1966} provides the following:
\begin{description}[before=\let\makelabel\descriptionlabel]
	\item[\textbf{(CH2)}\refstepcounter{desccount}\label{DESC: CH2}] Any closed subset $\mathcal{F}$ of $\mathcal{X}$ is a \textit{taut} subspace:
	\begin{equation}
		\label{EQ: TautnessOfCechCohomologies}
		\check{H}^{*}(\mathcal{F};G) = \varinjlim \check{H}^{*}(\mathcal{U};G),
	\end{equation}
	where the direct limit is taken with respect to the set of all open neighborhoods $\mathcal{U}$ of $\mathcal{F}$ in $\mathcal{X}$ directed downward by inclusion.
\end{description}

Next, we will need \v{C}ech cohomology module $\check{H}^{*}(\mathcal{A},\Gamma;G)$ for a pair $(\mathcal{A},\Gamma)$ of compact metric spaces $\Gamma \subset \mathcal{A}$. In this context, \cite[Theorem 7.6, Chapter IX]{EilenbergSteenrod2015} gives the following:
\begin{description}[before=\let\makelabel\descriptionlabel]
	\item[\textbf{(CH3)}\refstepcounter{desccount}\label{DESC: CH3}] There is a long exact sequence in cohomology for any pair of compact spaces.
\end{description}
Next, \cite[Lemma 3.1]{Morita1975} yields the following:
\begin{description}[before=\let\makelabel\descriptionlabel]
	\item[\textbf{(CH4)}\refstepcounter{desccount}\label{DESC: CH4}] If $\operatorname{dim}\mathcal{A} \leq m$, then $\check{H}^{m+k}(\mathcal{A};\Gamma;G) = 0$ for any $k=1,2,\ldots$.
\end{description}

The last two properties are used only once in the following lemma.
\begin{lemma}
	\label{LEM: BendixsonCohomL1}
	Let $\Gamma \subset \mathcal{A}$ be compact subsets of $\mathcal{X}$. If $\dim_{\operatorname{T}}\mathcal{A} \leq m$, then the inclusion $i_{\mathcal{A}} \colon \Gamma \to \mathcal{A}$ induces a surjective mapping $i^{*}_{\mathcal{A}} \colon \check{H}^{m}(\mathcal{A};G) \to \check{H}^{m}(\Gamma;G)$.
\end{lemma}
\begin{proof}
	By \nameref{DESC: CH3}, we have the exact sequence in cohomology of $(\mathcal{A},\Gamma)$:
	\begin{equation}
		\label{EQ: CechLongExactSequencePair}
		\cdots \rightarrow \check{H}^{m}(\mathcal{A};G) \overset{i^{*}_{\mathcal{A}}}{\rightarrow} \check{H}^{m}(\Gamma;G) \rightarrow \check{H}^{m+1}(\mathcal{A},\Gamma;G) \rightarrow \cdots.
	\end{equation}
	By \nameref{DESC: CH4}, we have that $\check{H}^{m+1}(\mathcal{A},\Gamma;G)=0$. Consequently, by exactness, $i^{*}_{\mathcal{A}}$ is surjective.
\end{proof}

For an open neighborhood $\mathcal{U}$ of $\mathcal{A}$, let $a_{\mathcal{U}} \colon \mathcal{A} \to \mathcal{U}$ be the inclusion mapping. The following abstract lemma contains a key argument common to all four cases considered below.
\begin{lemma}
	\label{LEM: BendixsonAbstractLem2}
	Let $\Gamma$ and $\mathcal{A}$ be compact subsets of $\mathcal{X}$. Consider a finitely generated subgroup $W$ of $\check{H}^{m}(\mathcal{A};G)$. Suppose there are sequences of metric spaces $\mathcal{B}_{k}$ and continuous mappings $T_{k} \colon \mathcal{B}_{k} \to \mathcal{X}$ and $j_{k} \colon \Gamma \to \mathcal{B}_{k}$ such that for any open neighborhood $\mathcal{U}$ of $\mathcal{A}$ there is a nonempty set $K(\mathcal{U})$ of $k=1,2,\ldots$ satisfying the following:
	\begin{enumerate}
		\item[1)] $T_{k} (\mathcal{B}_{k}) \subset \mathcal{U}$;
		\item[2)] $A_{k} \coloneq T_{k} \circ j_{k}$ maps $\Gamma$ into $\mathcal{A}$, so $A_{k} \colon \Gamma \to \mathcal{A}$.
	\end{enumerate}
	Then there exists a neighborhood $\mathcal{U}_{W}$ of $\mathcal{A}$ such that for any $k \in K(\mathcal{U}_{W})$ the image of $\check{H}^{m}(\mathcal{B}_{k};G)$ in $\check{H}^{m}(\Gamma;G)$ under $j^{*}_{k}$ contains $A^{*}_{k}(W)$.
\end{lemma}
\begin{proof}
	Let $v_{1},\ldots,v_{n}$ be a basis of $W$. By \nameref{DESC: CH2}, for any $i \in \{1,\ldots,n\}$ there exists an open neighborhood $\mathcal{U}_{i}$ of $\mathcal{A}$ such that $v_{i} \in a^{*}_{\mathcal{U}_{i}}(\check{H}^{m}(\mathcal{U}_{i};G))$. Let $\mathcal{U}_{W}$ be the intersection of all $\mathcal{U}_{i}$. Then $W \subset H_{W} \coloneq a^{*}_{\mathcal{U}_{W}}(\check{H}^{m}(\mathcal{U}_{W};G))$.
	
	Let $k \in K(\mathcal{U}_{W})$. By item 1), $T_{k}$ determines a mapping $T_{k, W} \colon \mathcal{B}_{k} \to \mathcal{U}_{W}$. Note that $T_{k, W} \circ j_{k} = a_{\mathcal{U}_{W}} \circ A_{k}$. Thus,
	\begin{equation}
		j^{*}_{k}(\check{H}^{m}(\mathcal{B}_{k};G)) \supset j^{*}_{k}(T^{*}_{k,W}(\check{H}^{m}(\mathcal{U}_{W};G))) = A^{*}_{k}(H_{W}) \supset A^{*}_{k}(W).
	\end{equation} 
\end{proof}

First, we illustrate the abstract argument by deriving the Bendixson criterion for attractors of semiflows. Note that it requires neither the injectivity of semiflow mappings nor the invariance of $\Gamma$.
\begin{theorem}
	\label{TH: BendixsonSemiflow}
	Let $\varphi$ be a continuous semiflow on $\mathcal{X}$, and let $\Gamma \subset \mathcal{A}$ be compact subsets of $\mathcal{X}$ such that for a positive integer $m$ we have the following:
	\begin{enumerate}
		\item[1)] $\mathcal{A}$ attracts a subset $\mathcal{B} \supset \Gamma$ of $\mathcal{X}$ with $\check{H}^{m}(\mathcal{B};G) = 0$;
		\item[2)] $\varphi^{t}(\Gamma) \subset \mathcal{A}$ for any $t \geq 0$;
		\item[3)] $\dim_{\operatorname{T}}\mathcal{A} \leq m$ (if $\Gamma \not= \mathcal{A}$).
	\end{enumerate}
	Then $\check{H}^{m}(\Gamma;G) = 0$.
\end{theorem}
\begin{proof}
	We consider Lemma \ref{LEM: BendixsonAbstractLem2} with $\mathcal{B}_{k} \coloneq \mathcal{B}$, $T_{k} \coloneq \restr{\varphi^{k}}{\mathcal{B}}$, and $j_{k}=i_{\mathcal{B}}$ given by the inclusion of $\Gamma$ into $\mathcal{B}$, which is correct by item 1). Note that $A_{k} \colon \Gamma \to \mathcal{A}$ is homotopic to $i_{\mathcal{A}}$ through $[0,1] \ni s \mapsto \restr{\varphi^{sk}}{\Gamma}$, which is well defined by item 2). By \nameref{DESC: CH1}, $A^{*}_{k} = i^{*}_{\mathcal{A}}$. By Lemma \ref{LEM: BendixsonAbstractLem2}, for any $w \in \check{H}^{m}(\mathcal{A};G)$ we have that
	\begin{equation}
		A^{*}_{k}(w) = i^{*}_{\mathcal{A}}(w) \in i^{*}_{\mathcal{B}}(\check{H}^{m}(\mathcal{B};G)) = \{0\}.
	\end{equation}
	Finally, item 3) and Lemma \ref{LEM: BendixsonCohomL1} ensure that $\check{H}^{m}(\Gamma;G) = i^{*}_{\mathcal{A}}(\check{H}^{m}(\mathcal{A};G)) = 0$.
\end{proof}

Next, we consider continuous-time extensions. For definitions, see Section \ref{SEC: Preliminaries}.
\begin{theorem}
	\label{TH: BendixsonCocycleContTime}
	Let $\wp \colon (\mathcal{Z},\pi) \to (\mathcal{Q}, \vartheta)$ be an extension of a continuous flow $(\mathcal{Q}, \vartheta)$ to a continuous semiflow $(\mathcal{Z},\pi)$. Fix $q \in \mathcal{Q}$, and let $\Gamma \subset \mathcal{A} \subset \mathcal{Z}_{q}$ be compact subsets in the fiber $\mathcal{Z}_{q}$ over $q$ such that for a positive integer $m$ we have the following:
	\begin{enumerate}
		\item[1)] $\mathcal{A}$ pullback attracts over $q$ a family of subsets $\{\mathcal{B}_{s} \}_{s \leq 0}$ such that $\mathcal{B}_{0} \supset \Gamma$, $\mathcal{B}_{s} \subset \mathcal{Z}_{\vartheta^{s}(q)}$, and $\check{H}^{m}(\mathcal{B}_{s};G) = 0$ for any $s \leq 0$;
		\item[2)] there exist continuous mappings $j_{s} \colon \Gamma \to \mathcal{B}_{s}$, where $s \leq 0$, such that $j_{0} = i_{\mathcal{B}_{0}} \colon \Gamma \to \mathcal{B}_{0}$ is the inclusion map, $\hat{A}_{-t} \coloneq \pi^{t} \circ j_{-t}$ maps $\Gamma$ into $\mathcal{A}$ for any $t \geq 0$, and $(-\infty,0] \times \Gamma \ni (s,x) \mapsto \hat{A}_{s}(x) \in \mathcal{A}$ is continuous;
		\item[3)] $\dim_{\operatorname{T}}\mathcal{A} \leq m$ (if $\Gamma \not= \mathcal{A}$).
	\end{enumerate}
	Then $\check{H}^{m}(\Gamma;G) = 0$.
\end{theorem}
\begin{proof}
	We consider Lemma \ref{LEM: BendixsonAbstractLem2} with $\mathcal{X} \coloneq \mathcal{Z}_{q}$, $\mathcal{B}_{k} \coloneq \mathcal{B}_{-k}$, $T_{k} \coloneq \restr{\pi^{k}}{\mathcal{B}_{-k}}$, and $j_{k} \coloneq j_{-k}$. Note that $A_{k} \coloneq \hat{A}_{-k} \colon \Gamma \to \mathcal{A}$ is homotopic to $i_{\mathcal{A}}$ through $[0,1] \ni s \mapsto \hat{A}_{-sk}$, which is well defined by item 2). Thus, $A^{*}_{k} = i^{*}_{\mathcal{A}}$ and we can proceed as in the proof of Theorem \ref{TH: BendixsonSemiflow}.
\end{proof}
\begin{remark}
	\label{REM: BendixsonForContTimeCocyclesRemark}
	In the context of the above theorem, if there are compact subsets $\Gamma_{s} \subset \mathcal{B}_{s}$ such that $\pi^{t}$ takes $\Gamma_{-t}$ homeomorphically onto $\Gamma$, then $j_{-t}$ can be defined by the inverse mapping considered with values in $\mathcal{B}_{-t}$. In this case, $\hat{A}_{-t} = i_{\mathcal{A}}$, and hence the assumptions of item 2) are trivially satisfied.
\end{remark}

Now we treat discrete-time systems, where the assumptions on the invariance of $\Gamma$ and injectivity of mappings on $\Gamma$ are necessary. We start with the case of a single mapping.
\begin{theorem}
	\label{TH: BendixsonMapping}
	Let $T$ be a continuous mapping of $\mathcal{X}$ into itself, and let $\Gamma \subset \mathcal{A}$ be compact subsets of $\mathcal{X}$ such that for a positive integer $m$ we have the following:
	\begin{enumerate}
		\item[1)] $\mathcal{A}$ attracts a subset $\mathcal{B} \supset \Gamma$ of $\mathcal{X}$ with $\check{H}^{m}(\mathcal{B};G) = 0$;
		\item[2)] $T$ takes $\Gamma$ homeomorphically into itself, or at least $T(\Gamma) \subset \Gamma$ and $T^{*}$ is injective on $\check{H}^{m}(\Gamma;G)$;
		\item[3)] $\dim_{\operatorname{T}}\mathcal{A} \leq m$ (if $\Gamma \not= \mathcal{A}$);
	\end{enumerate}
	Then $\check{H}^{m}(\Gamma;G)=0$.
\end{theorem}
\begin{proof}
	We consider Lemma \ref{LEM: BendixsonAbstractLem2} with $\mathcal{B}_{k} \coloneq \mathcal{B}$, $T_{k} \coloneq \restr{T^{k}}{\mathcal{B}}$, and $j_{k}=i_{\mathcal{B}}$ given by the inclusion of $\Gamma$ into $\mathcal{B}$, which is correct by item 1). Note that $F_{k} \coloneq \restr{T^{k}}{\Gamma}$ is the mapping $T_{k} \circ j_{k}$ considered as the self-map of $\Gamma$, which is correct by item 2). Then $A_{k} = i_{\mathcal{A}} \circ F_{k}$, where $i_{\mathcal{A}}$ is the inclusion of $\Gamma$ into $\mathcal{A}$.
	
	Suppose that $v \in \check{H}^{m}(\Gamma;G)$. By item 3) and Lemma \ref{LEM: BendixsonCohomL1}, there exists $w \in \check{H}^{m}(\mathcal{A};G)$ such that $v = i^{*}_{\mathcal{A}}(w)$. Then Lemma \ref{LEM: BendixsonAbstractLem2} guarantees that there exists $k$ such that
	\begin{equation}
		A^{*}_{k}(w) = F^{*}_{k}(i^{*}_{\mathcal{A}}(w)) = F^{*}_{k}(v) \in i^{*}_{\mathcal{B}}(\check{H}^{m}(\mathcal{B};G)) = \{0\}.
	\end{equation}
	Since $F^{*}_{k}$ is injective by item 2), we must have $v = 0$.
\end{proof}
\begin{remark}
	The assumption on injectivity in item 3) cannot be discarded. A counterexample is constructed in \cite[p.~465]{LiMuldowney1995LowBounds}, where $\Gamma$ is a topological circle.
\end{remark}

Finally, we treat discrete-time extensions, which without loss of generality can be identified with sequences of mappings.
\begin{theorem}
	\label{TH: BendixsonCocyclesDiscreteTime}
	For $k=1,2,\ldots$, let $T_{k}$ be a continuous mapping of a metric space $\mathcal{B}_{k}$ into $\mathcal{X}$, and let $\Gamma \subset \mathcal{A}$ be compact subsets of $\mathcal{X}$ such that for a positive integer $m$ we have the following:
	\begin{enumerate}
		\item[1)] For any neighborhood $\mathcal{U}$ of $\mathcal{A}$, there is $k$ such that $T_{k}(\mathcal{B}_{k}) \subset \mathcal{U}$ and $\check{H}^{m}(\mathcal{B}_{k};G) = 0$;
		\item[2)] For any $k$, there is a continuous mapping $j_{k} \colon \Gamma \to \mathcal{B}_{k}$ such that $A_{k} \coloneq T_{k} \circ j_{k}$ takes $\Gamma$ into itself homeomorphically, or at least such that $A_{k}(\Gamma) \subset \Gamma$ and $A^{*}_{k}$ is injective on $\check{H}^{m}(\Gamma;G)$.
		\item[3)] $\dim_{\operatorname{T}}\mathcal{A} \leq m$ (if $\Gamma \not= \mathcal{A}$);
	\end{enumerate}
	Then $\check{H}^{m}(\Gamma;G)=0$.
\end{theorem}
\begin{proof}
	We consider Lemma \ref{LEM: BendixsonAbstractLem2} with $\mathcal{B}_{k}$, $T_{k}$, and $j_{k}$ as in the statement. Let $F_{k}$ be $A_{k}$ considered as the self-mapping of $\Gamma$, which is correct by item 2). Clearly, $A_{k} = i_{\mathcal{A}} \circ F_{k}$. Reasoning as in the proof of Theorem \ref{TH: BendixsonMapping}, for any $v \in \check{H}^{m}(\Gamma;G)$ we obtain $w \in \check{H}^{m}(\mathcal{A};G)$ such that $v=i^{*}_{\mathcal{A}}(w)$ and $A^{*}_{k}(w) = F^{*}_{k}(v) = 0$. By item 2), $F^{*}_{k}$ is injective, and hence $v=0$.
\end{proof}
\begin{remark}
	Similarly to Remark \ref{REM: BendixsonForContTimeCocyclesRemark}, if there are compact subsets $\Gamma_{k}$ of $\mathcal{B}_{k}$ such that $T_{k}$ takes $\Gamma_{k}$ homeomorphically onto $\Gamma$, we can define $j_{k}$ using the inverse mapping, and the assumption of item 2) is satisfied.
\end{remark}

When $\Gamma$ is homeomorphic to the boundary of an open bounded subset of $\mathbb{R}^{m+1}$, then the group $\check{H}^{m}(\Gamma;\mathbb{Z})$ is not trivial by Alexander duality. This case is treated in \cite{LiMuldowney1995LowBounds} with additional assumptions on $\Gamma$ concerning its local geometry in the underlying Banach space. Theorem \ref{TH: BendixsonMapping} shows that no additional assumptions are actually needed. Moreover, the semiflow case in \cite{LiMuldowney1995LowBounds} reduces to maps. This leads to another redundant assumption that the semiflow maps are injective or that $\Gamma$ is invariant. Theorem \ref{TH: BendixsonSemiflow} relaxes these assumptions as well. 

\begin{remark}
	Let us also mention that on \cite[p.~458]{LiMuldowney1995LowBounds} it is misleadingly claimed that the proof of \cite[Theorem 5]{Smith1986HD} ``tacitly assumes'' injectivity of the mapping on $\Gamma$. Although the injectivity condition does indeed make the statement true, the proof actually fails in the last paragraph, where an attempt is made to apply a homotopy argument that goes beyond the theorem's formal assumptions. Thus, the right ``tacit assumption'' would be the requirement that the mapping on the attractor is homotopic to the identity. This is precisely what we have for semiflows, and the interest of \cite{Smith1986HD} in applications to ODEs only confirms this.
\end{remark}

Another interesting case involving nontrivial cohomology arises when $\Gamma$ is a compact abelian group of topological dimension $m$. In this case, $\check{H}^{m}(\Gamma;\mathbb{Z})$ is always nontrivial, see \cite[Chapter 5]{HofmannMostert1973}. This will be used in Section \ref{SEC: AlmostPeriodicBendixsonRefinement} to refine certain results for almost periodic flows.

Although the Bendixson criterion for attractors is formulated as a negative criterion, it can be used to obtain positive results. Specifically, the paper \cite{LiMuldowney1995LowBounds} investigates lower bounds on the dimension of $\mathcal{A}$. Namely, if the conditions of the theorems stated above are satisfied, except for item 3), and $\check{H}^{m}(\Gamma;G) \not=0$, then we necessarily have that $\dim_{\operatorname{T}}\mathcal{A} \geq m+1$.

We present another positive result in the following theorem.
\begin{theorem}
	\label{TH: BendixsonCriterionSemiflowsPositiveLocalizationResult}
	Let $\varphi$ be a continuous semiflow in $\mathcal{X}$. Suppose that $\Gamma$ is a compact orbitally stable\footnote{That is, for any neighborhood $\mathcal{U}$ of $\Gamma$ there exists another neighborhood $\mathcal{V}$ such that $\varphi^{t}(\mathcal{V})$ lies in $\mathcal{U}$ for all $t \geq 0$.} invariant set with $\check{H}^{m}(\Gamma;G) \not= 0$ for some $m$, and there is a compact set $\Omega \supset \Gamma$ with $\check{H}^{m}(\Omega;G) = 0$. Then there exists a point $x \in \Omega$ that is not attracted by $\Gamma$.
\end{theorem}
\begin{proof}
	Indeed, if any point of $\Omega$ is attracted by $\Gamma$, then, by virtue of continuity, compactness, and orbital stability, $\Gamma$ attracts $\Omega$. According to Theorem \ref{TH: BendixsonSemiflow} with $\mathcal{A}=\Gamma$ and $\mathcal{B} = \Omega$, we must have $\check{H}^{m}(\Gamma;G) = 0$ in this case, which leads to a contradiction.
\end{proof}
\begin{remark}
	A similar statement holds for continuous mappings $T$ under the additional condition that $T$ maps $\Gamma$ homeomorphically onto itself.
\end{remark}

In particular, we have the following corollary.
\begin{corollary}
	Let $\varphi$ be a continuous semiflow in a normed vector space $\mathcal{X}$, and let $\Gamma$ be an orbitally stable periodic orbit of $\varphi$. Then for any $x \in \Gamma$ there exists $y \in \Gamma$ such that the line segment between $x$ and $y$ contains a point, which is not attracted by $\Gamma$.
\end{corollary}
\begin{proof}
	Fix $x \in \Gamma$, and let $\Omega$ be formed by all the line segments from $x$ to other points of $\Gamma$. For $z \in \Omega$ and $s \in [0,1]$, define $H(z, s) \coloneq (1-s)z + s x$. It is not hard to see that $H(\Omega,s) \subset \Omega$ for any $s \in [0,1]$, and $H$ contracts $\Omega$ to the one-point set $\{x\}$. Thus, $\check{H}^{1}(\Omega;\mathbb{Z}) = 0$, whereas $\check{H}^{1}(\Gamma;\mathbb{Z}) = \mathbb{Z}$. Applying Theorem \ref{TH: BendixsonCriterionSemiflowsPositiveLocalizationResult} yields the desired result.
\end{proof}

The above result is of interest from the perspective of numerical analysis. Specifically, if a periodic attractor $\Gamma$ has been found, one can track trajectories starting on the two-dimensional set $\Omega$ to reveal other nontrivial behavior in the system. This does not necessarily lead to the discovery of other attractors. For example, it may turn out that the only point not attracted to $\Gamma$ lies on the stable manifold of an unstable equilibrium. Numerically, we can observe only transient processes in the vicinity of this equilibrium, which may lead to a return to $\Gamma$. This is precisely what occurs in the case of an Andronov--Hopf bifurcation.
\section{Semiflows with one-dimensional attractors}
\label{SEC: GlobalStabilitySemiflowsBendixson}
In this section, we aim to prove Theorem \ref{TH: BendixsonGlobalStabilityOfSemiflows}. For what follows, let $(\mathcal{X}, \varphi)$ be a continuous semiflow in a metric space $\mathcal{X}$ with distance function $d_{\mathcal{X}}$.

Recall that a point $x \in \mathcal{X}$ is called \textit{nonwandering} if there exist a sequence of points $x_{1},x_{2},\ldots \in \mathcal{X}$ and a sequence of arbitrarily large times $t_{1},t_{2}, \ldots$ such that $x_{k} \to x$ and $\varphi^{t_{k}}(x_{k}) \to x$ as $k \to \infty$.

A point $x \in \mathcal{X}$ is called \textit{chain recurrent} if for any $\varepsilon>0$ and $T > 0$ there exist times $s_{1}, \ldots, s_{n} \geq T$ and points $x_{0}, \ldots, x_{n} \in \mathcal{K}$ with $x_{0} = x_{n} = x$ such that
\begin{equation}
	d_{\mathcal{X}}(\varphi^{s_{j}}(x_{j-1}), x_{j}) < \varepsilon \qquad \text{for all} \ j \in \{ 1,\ldots,n\}.
\end{equation}
In this context, the sequence $x_{0},\ldots, x_{n}$ is called an $(\varepsilon,T)$-\textit{chain} with \textit{jump times} $s_{1},\ldots, s_{n}$.

Note that both concepts describe a property not of the trajectory of the point $x$ itself, but rather of the trajectories of points within small neighborhoods of $x$. Because of this, a point $x$ may be chain recurrent (resp. nonwandering) with respect to $(\mathcal{X},\varphi)$, but lose this property when restricting $\varphi$ to an invariant subset of $\mathcal{X}$ containing $x$. 

This can be illustrated by means of $\omega$-limit sets. Recall that the $\omega$-limit set $\omega(x)$ of $x \in \mathcal{X}$ is the set of all $y \in \mathcal{Y}$ such that $\varphi^{t_{k}}(x) \to y$ as $k \to \infty$ for some sequence $t_{1},t_{2},\ldots$ of arbitrary large times. Assuming that the positive semiorbit of $x$ has compact closure, it is well known that $\omega(x)$ is a nonempty compact connected invariant set that attracts $x$. Moreover, any $y \in \omega(x)$ is nonwandering with respect $(\mathcal{X},\varphi)$ (directly by definition) and chain recurrent with respect to $(\omega(x),\varphi)$ (see Lemma \ref{LEM: NonwanderingAreChainRecurrent}). When there is a complete trajectory $\gamma$ of $\varphi$ passing through $x$ (i.e., $\gamma(0) = x$) such that $\gamma((-\infty,0])$ has compact closure, similar conclusions can be deduced for the $\alpha$-limit set $\alpha(\gamma)$ of $x$ with respect to $\gamma$.
\begin{remark}
	In terms of complete trajectories, we also use the notation $\omega(\gamma) \coloneq \omega(\gamma(0))$. When there is a unique complete trajectory passing through a given point $x$ (e.g., when $\varphi$ is a flow), then we also write $\alpha(x)$ instead of $\alpha(\gamma)$.
\end{remark}
 
\subsection{From semiflows to flows}

Let us begin with the following lemma, which illustrates how external nonwandering behavior translates into internal chain recurrence in the presence of attraction. In this form, the lemma can be applied to nonwandering points belonging to an attractor of compact sets (in particular, to any point in an $\omega$-limit set), as well as to any point in an $\alpha$-limit set.
\begin{lemma}
	\label{LEM: NonwanderingAreChainRecurrent}
	Let $\varphi$ be a semiflow in a metric space $\mathcal{X}$, and let $\mathcal{A}$ be a compact invariant set. Suppose that for some $q \in \mathcal{A}$ there exist a sequence of points $p_{1},p_{2},\ldots \in \mathcal{X}$ and a sequence of arbitrary large times $t_{1},t_{2},\ldots$ such that $p_{k} \to q$ and $\varphi^{t_{k}}(p_{k}) \to q$ as $k \to \infty$, with the following property being satisfied:
	\begin{equation}
		\label{EQ: NonwanderingAttractingAreChainRecurrent}
		\lim_{k \to \infty}\sup_{s \in [0,t_{k}]} d_{\mathcal{X}}(\varphi^{s}(p_{k}), \mathcal{A}) = 0.
	\end{equation}
	Then $q$ is chain recurrent with respect to $(\mathcal{A},\varphi)$.
\end{lemma}
\begin{proof}
	Let $T > 0$ and $\varepsilon > 0$ be fixed. By continuity of $\varphi$ and compactness of $\mathcal{A}$, there exist $\delta>0$ and a finite number of points $q'_{1},\ldots, q'_{m} \in \mathcal{A}$ such that for any $p \in \mathcal{O}_{\delta}(\mathcal{A})$ there exists $i \in \{1,\ldots,m\}$ such that
	\begin{equation}
		\label{EQ: NonwanderingGivesChainRec1}
		d_{\mathcal{X}}(\varphi^{t}(p), \varphi^{t}(q'_{i})) < \varepsilon/2 \qquad \text{for any} \ t \in [0,2T].
	\end{equation}
	
	Let $p_{1},p_{2},\ldots$ and $t_{1},t_{2},\ldots$ be as in the statement. By \eqref{EQ: NonwanderingAttractingAreChainRecurrent}, there exists $N$ such that
	\begin{equation}
		d_{\mathcal{X}}( \varphi^{t}(p_{k}), \mathcal{A}) < \delta \qquad \text{for all} \ k \geq N \ \text{and} \ t \in [0, t_{k}].
	\end{equation}
	
	Fix $k \geq N$ such that $d_{\mathcal{X}}(\varphi^{T}(p_{k}), \varphi^{T}(q)) < \delta$, $d_{\mathcal{X}}(\varphi^{t_{k}}(p_{k}), q) < \varepsilon/2$, and $t_{k} > 2T$. Set $s_{1} \coloneq T$. Then there exist an integer $n \geq 2$ and times $s_{2},\ldots, s_{n}$ in $[T,2T]$ such that $t_{k} = s_{1} + \cdots + s_{n}$. Set $\Sigma_{j} \coloneq s_{1} + \cdots + s_{j}$, By \eqref{EQ: NonwanderingGivesChainRec1}, for any $j \in \{1,\ldots, n-1\}$ there exists $q_{j} \in \{ q'_{1}, \ldots, q'_{m} \}$ such that
	\begin{equation}
		d_{\mathcal{X}}(\varphi^{\Sigma_{j} + t}(p_{k}), \varphi^{t}(q_{j})) < \varepsilon/2 \qquad \text{for any} \ t \in [0,s_{j+1}].
	\end{equation}
	Using this with $t=0$ or $t=s_{j+1}$ and the choice of $k$, we obtain that
	\begin{equation}
		d_{\mathcal{X}}(\varphi^{s_{j}}(q_{j-1}), q_{j}) \leq d_{\mathcal{X}}(\varphi^{\Sigma_{j}}(p_{k}), \varphi^{s_{j}}(q_{j})) +  d_{\mathcal{X}}(\varphi^{\Sigma_{j}}(p_{k}), q_{j}) < \varepsilon
	\end{equation}
	for any $j \in \{1,\ldots,n\}$, where $q_{0}=q_{n} = q$. Since $T$ and $\varepsilon$ were arbitrary, this establishes that $q$ is chain recurrent with respect to $(\mathcal{A},\varphi)$.
\end{proof}

Now we consider a semiflow $(\mathcal{K},\varphi)$ on a compact metric space $\mathcal{K}$ with surjective time-$t$ mappings. As discussed in Section \ref{SEC: Preliminaries}, there is a flow $(\hat{K},\hat{\varphi})$ called the natural extension of $(\mathcal{K},\varphi)$. We will need the following relation between chain recurrent points.
\begin{lemma}\cite[Theorem 3.4]{HeZhouLv2002}
	\label{LEM: ChainRecurrenceLiftInverseLimit}
	Any chain recurrent point of $(\mathcal{K},\varphi)$ lifts to a chain recurrent point of $(\hat{\mathcal{K}},\hat{\varphi})$ and vice versa.
\end{lemma}

The following lemma is the final piece justifying the transition from $(\mathcal{K},\varphi)$ to its natural extension.
\begin{lemma}
	\label{LEM: DimAndCohomInverseLimit}
	For $(\hat{\mathcal{K}},\hat{\varphi})$ as above, we have the following:
	\begin{enumerate}
		\item[(a)] $\dim_{\operatorname{T}}\hat{\mathcal{K}} \leq \dim_{\operatorname{T}} \mathcal{K}$;
		\item[(b)] If $\check{H}^{1}(\mathcal{K};\mathbb{Z})=0$, then $\check{H}^{1}(\hat{\mathcal{K}};\mathbb{Z}) = 0$.
	\end{enumerate}
\end{lemma}
\begin{proof}
	This follows from the fact that $\hat{\mathcal{K}}$ is the inverse limit of the sequence of copies $\mathcal{K}_{i}$ of $\mathcal{K}$ with bonding maps $f_{i+1} \coloneq \varphi^{1} \colon \mathcal{K}_{i+1} \to \mathcal{K}_{i}$ and the behavior of dimension and cohomology under such limits. For (a), see \cite[Theorem 1.13.4]{Engelking1978}. For (b), see \cite[Chapter X, Theorem 3.1]{EilenbergSteenrod2015}.
\end{proof}

\subsection{Chain recurrent points of flows on compact spaces}

In this subsection, we fix a continuous flow $\varphi$ on a compact metric space $\mathcal{K}$ with distance function $d_{\mathcal{K}}$. Let $\operatorname{Eq}(\mathcal{K},\varphi)$ denote the set of equilibrium points of $(\mathcal{K},\varphi)$, i.e., $q \in \mathcal{K}$ such that $\varphi^{t}(q)=q$ for any $t \in \mathbb{R}$.

The following lemma is the well-known Whitney--Bebutov theorem on local sections, see, e.g., \cite[Lemma 3.2]{GutmanJinTsukamoto2019}.
\begin{lemma}
	\label{LEM: WhitneyBebutovSections}
	Let $q \notin \operatorname{Eq}(\mathcal{K},\varphi)$. Then there exist a compact subset $\Sigma \subset \mathcal{K}$ containing $q$ and $\delta > 0$ such that
	\begin{equation}
		\label{EQ: WhitneyBebutovLocalSectionMapping}
		\Sigma \times [-\delta,\delta] \ni (p,t) \mapsto \varphi^{t}(p) \in \mathcal{K}
	\end{equation}
	is an injective continuous mapping whose image contains an open neighborhood of $q$.
\end{lemma}

The following lemma refines the structure of local sections for the case of one-dimensional phase spaces.
\begin{lemma}
	\label{LEM: WhitneyBebutovSectionZeroDim}
	Suppose that $\dim_{\operatorname{T}}\mathcal{K} = 1$ and $q \notin \operatorname{Eq}(\mathcal{K},\varphi)$. Then there exist a compact zero-dimensional set $\Sigma_{0}$, $\delta_{0} > 0$, and a homeomorphism $F \colon \Sigma_{0} \times [-\delta_{0}, \delta_{0}] \to \mathcal{N}_{0}$ onto a compact subset $\mathcal{N}_{0}$ of $\mathcal{K}$ such that
	\begin{enumerate}
		\item[1)] $F(p,t) = \varphi^{t}(p)$ for $p \in \Sigma_{0}$ and $t \in [-\delta_{0},\delta_{0}]$, and $F(p_{0},0) = q$ for some $p_{0} \in \Sigma_{0}$;
		\item[2)] $F(\Sigma_{0} \times (-\delta_{0},\delta_{0})) \eqcolon \mathcal{U}_{0}$ is an open neighborhood of $q$ in $\mathcal{K}$;
		\item[3)] $\partial \mathcal{U}_{0} = F(\Sigma_{0} \times \{ \pm \delta_{0} \})$;
		\item[4)] $\mathcal{N}_{0} \cap \operatorname{Eq}(\mathcal{K},\varphi) = \emptyset$. 
	\end{enumerate}
\end{lemma}
\begin{proof}
	Let $\Sigma$ and $\delta$ be as in Lemma \ref{LEM: WhitneyBebutovSections}, and let $F$ denote the mapping in \eqref{EQ: WhitneyBebutovLocalSectionMapping}. Set $\mathcal{N} \coloneq F(\Sigma \times [-\delta, \delta])$. Since $\mathcal{N}$ is one-dimensional, $\Sigma$ must be zero dimensional. This follows from the Hurewicz theorem: $\dim_{\operatorname{T}}(\Sigma \times [-\delta,\delta]) = \dim_{\operatorname{T}}\Sigma + 1$, see \cite[Problem 1.9.E (b)]{Engelking1978}. In particular, $\Sigma$ has a base of clopen subsets.
	
	Let $\mathcal{U}$ be the interior of $\mathcal{N}$ in $\mathcal{K}$. Then $q \in \mathcal{U}$, $q = F(p_{0},0)$ for some $p_{0} \in \Sigma$, and $\mathcal{W} \coloneq F^{-1}(\mathcal{U})$ is an open subset in $\Sigma \times [-\delta, \delta]$ containing $(p_{0},0)$.	Let $\Sigma_{0}$ be any clopen subset of $\Sigma$ containing $p_{0}$ and $\delta_{0} \in (0, \delta)$ such that $\Sigma_{0} \times [-\delta_{0}, \delta_{0}] \subset \mathcal{W}$. Clearly, item 1) is satisfied for such defined objects. Then item 4) immediately follows from this.
	
	Note that $\mathcal{N}_{0} \coloneq F(\Sigma_{0} \times [-\delta_{0},\delta_{0}])$ is compact and lies in $\mathcal{U}$. Since $\Sigma_{0} \times (-\delta_{0},\delta_{0})$ is open in $\Sigma \times [-\delta,\delta]$, its image $\mathcal{U}_{0}$ under $F$ is open in $\mathcal{N}$. Since $\mathcal{U}_{0} \subset \mathcal{U} \subset \mathcal{N}$, it is also open in $\mathcal{U}$ and hence in $\mathcal{K}$. This shows item 2).
	
	Clearly, the closure of $\mathcal{U}_{0}$ is $\mathcal{N}_{0}$. Consequently, $\partial \mathcal{U}_{0} = \mathcal{N}_{0} \setminus \mathcal{U}_{0} = F(\Sigma_{0} \times \{\pm \delta_{0} \})$. This shows item 3).
\end{proof}

The following statement is well-known (in fact, the converse to it also holds), but we present its elementary proof for the sake of completeness.
\begin{lemma}
	\label{LEM: CohomologyZeroLiftingFtoG}
	Let $f \colon \mathcal{K} \to \mathbb{T}^{1}$ be a continuous function from $\mathcal{K}$ to the unit circle $\mathbb{T}^{1} \subset \mathbb{C}$. If $\check{H}^{1}(\mathcal{K};\mathbb{Z}) = 0$, then there exists a continuous function $g \colon \mathcal{K} \to \mathbb{R}$ such that $f(q) = e^{i2\pi g(q)}$ for any $q \in \mathcal{K}$.
\end{lemma}
\begin{proof}
	Let $\delta > 0$ be such that $|f(p)-f(q)| < 2$ provided that $p, q \in \mathcal{K}$ and $d_{\mathcal{K}}(p,q) < \delta$. Let $\mathcal{U}_{1}, \ldots, \mathcal{U}_{n}$ be a finite open cover of $\mathcal{K}$ by sets with diameters $\leq \delta$. For each $j \in \{1,\ldots,n\}$, there exists a continuous function $g_{j} \colon \mathcal{U}_{j} \to \mathbb{R}$ (a local lift) such that $f(q)=e^{i2\pi g_{j}(q)}$ for all $q \in \mathcal{U}_{j}$. On each intersection $\mathcal{U}_{k} \cap \mathcal{U}_{l}$, the difference $g_{k} - g_{l}$ is a locally constant integer-valued function, thereby defining a \v{C}ech $1$-cocycle with values in $\mathbb{Z}$. Since $\check{H}^{1}(\mathcal{K};\mathbb{Z}) = 0$, this cocycle is a coboundary, i.e., there exist functions $h_{j} \colon \mathcal{U}_{j} \to \mathbb{Z}$ such that $g_{k} - g_{l} = h_{k} - h_{l}$. Thus, the functions $g_{j}-h_{j}$ glue together into the desired $g$.
\end{proof}

Next, we construct local Lyapunov functions corresponding to local sections.
\begin{lemma}
	\label{LEM: ZeroDimFlowLocalGFunction}
	Under the hypotheses and notations of Lemma \ref{LEM: WhitneyBebutovSectionZeroDim}, let $\check{H}^{1}(\mathcal{K};\mathbb{Z}) = 0$. Then there exists a continuous function $g \colon \mathcal{K} \to \mathbb{R}$ such that
	\begin{enumerate}
		\item[1)] For any $p \in \mathcal{K}$, the function $t \mapsto g(\varphi^{t}(p))$ does not decrease;
		\item[2)] If $\varphi^{t_{0}}(p) \in \mathcal{U}_{0}$ for some $t_{0} \in \mathbb{R}$ and $p \in \mathcal{K}$, then the function $t \mapsto g(\varphi^{t}(p))$ strictly increases in a neighborhood of $t_{0}$. More precisely, it strictly increases on the connected component of $\{ t \in \mathbb{R} \colon \varphi^{t}(p) \in \mathcal{U}_{0} \}$ containing $t_{0}$, which has the form $(s,s+2\delta_{0})$ for some $s$, and $g(\varphi^{s+2\delta_{0}}(p)) - g(\varphi^{s}(p))= 1$;
		\item[3)] $g$ takes integer values on $\mathcal{K} \setminus \mathcal{U}_{0}$.
	\end{enumerate}
\end{lemma}
\begin{proof}
	Fix a continuous and strictly increasing function $\theta \colon [-\delta_{0}, \delta_{0}] \to [0,1]$ such that $\theta(-\delta_{0}) = 0$ and $\theta(\delta_{0}) = 1$. For $p \in \mathcal{K}$, define
	\begin{equation}
		\label{EQ: LocalLyapFunctionFromZeroDimSection}
		f(p) \coloneq \begin{cases}
			e^{i2\pi \theta(t)} \qquad &\text{if} \ p=F(p_{0},t) \ \text{with} \ (p_{0},t) \in \Sigma_{0} \times [-\delta_{0}, \delta_{0}],\\
			1 \qquad &\text{if} \ p \in \mathcal{K} \setminus \mathcal{U}_{0}.
		\end{cases}
	\end{equation}
	Each case determines a continuous function on one of the compact subsets $\mathcal{N}_{0}$ and $\mathcal{K} \setminus \mathcal{U}_{0}$. By item 3) of Lemma \ref{LEM: WhitneyBebutovSectionZeroDim}, these sets intersect at $\partial \mathcal{U}_{0}= F(\Sigma_{0} \times \{\pm \delta_{0} \})$, where both formulas yield the same result. Thus, $f$ is well-defined. By Lemma \ref{LEM: CohomologyZeroLiftingFtoG}, there exists $g \colon \mathcal{K} \to \mathbb{R}$ such that $f(p) = e^{i2\pi g(p)}$ for any $p \in \mathcal{K}$. Clearly, $g$ takes integer values on $\mathcal{K} \setminus \mathcal{U}_{0}$ and $g(F(p_{0},t)) = \theta(t) + m(p_{0})$ for $p_{0} \in \Sigma_{0}$ and $t \in [-\delta_{0},\delta_{0}]$, where $m$ takes integer values on $\Sigma_{0}$. In particular, item 3) is valid.
	
	Fix $p \in \mathcal{K}$ and let $\mathcal{I}$ be a connected component of $\{ t \in \mathbb{R} \colon \varphi^{t}(p) \in \mathcal{U}_{0} \}$. Since
	\begin{equation}
		\mathcal{I} \ni s \mapsto F^{-1}(\varphi^{s}(p)) \eqcolon (p_{0}(s), t(s)) \in \Sigma_{0} \times (-\delta_{0}, \delta_{0})
	\end{equation}
	is continuous, and $\Sigma_{0}$ is zero-dimensional, $p_{0}(s) \equiv p_{\mathcal{I}}$ on $\mathcal{I}$ for some $p_{\mathcal{I}} \in \Sigma_{0}$. Thus, $\mathcal{I} = (s, s + 2\delta_{0})$ and $p_{\mathcal{I}} = \varphi^{s+\delta_{0}}(p)$ for some $s$, and
	\begin{equation}
		g(\varphi^{t}(p)) = m(p_{\mathcal{I}}) + \theta(t - s - \delta_{0}) \quad \text{for} \ t \in \mathcal{I}.
	\end{equation}
	This establishes item 2).
	
	For item 1), choose two consecutive times $t_{1} < t_{2}$. Since the connected components $\mathcal{I}$ are disjoint open intervals of length $2\delta_{0}$, only at most a finite number of them may intersect the segment $[t_{1},t_{2}]$. In the complement, there is a finite number of disjoint closed intervals, on which the trajectory of $p$ belongs to $\mathcal{K} \setminus \mathcal{U}_{0}$. By item 3), $t \mapsto g(\varphi^{t}(p))$ is constant on each closed interval. Thus, as $t$ varies in $[t_{1},t_{2}]$, $g(\varphi^{t}(p))$ cannot decrease.
\end{proof}
\begin{remark}
	In Lemma \ref{LEM: ZeroDimFlowLocalGFunction}, the fact that $\Sigma_{0}$ is zero-dimensional is used to ensure the continuity of $f$ in \eqref{EQ: LocalLyapFunctionFromZeroDimSection}, and the key argument relies on the simple structure of the boundary $\partial \mathcal{U}_{0} = F(\Sigma_{0} \times \{\pm \delta_{0}\})$. In the case of general sections, continuity is lost precisely due to the appearance of other boundary components.
\end{remark}

In the next lemma, we glue local Lyapunov functions into a global Lyapunov function.
\begin{lemma}
	\label{LEM: ZeroDimFlowGlobalFunction}
	Under the hypotheses of Lemma \ref{LEM: ZeroDimFlowLocalGFunction}, there exists a continuous function $V \colon \mathcal{K} \to \mathbb{R}$ such that
	\begin{enumerate}
		\item[1)] For any $p \in \mathcal{K}$, $V(\varphi^{t}(p)) \geq V(p)$ for all $t \geq 0$, and $V(\varphi^{t}(p)) > V(p)$ provided that $t > 0$ and $p \notin \operatorname{Eq}(\mathcal{K},\varphi)$;
		\item[2)] $V(\operatorname{Eq}(\mathcal{K},\varphi))$ is a compact set of Lebesgue measure zero.
	\end{enumerate}
\end{lemma}
\begin{proof}
	For $q \notin \mathcal{K} \setminus \operatorname{Eq}(\mathcal{K},\varphi)$, Lemma \ref{LEM: WhitneyBebutovSectionZeroDim} gives an open flow-box $\mathcal{U}_{0}=\mathcal{U}_{0}(q)$ of $q$. Since the complement $\mathcal{K} \setminus \operatorname{Eq}(\mathcal{K},\varphi)$ is an open subset, it can be covered by a countable many such open flow-boxes $\mathcal{U}_{1},\mathcal{U}_{2},\ldots$. Let $g_{k}$ be the local Lyapunov function corresponding to the flow-box $\mathcal{U}_{k}$ by Lemma \ref{LEM: ZeroDimFlowLocalGFunction}.
	
	Since the set $\operatorname{Eq}(\mathcal{K},\varphi)$ is compact, each $g_{k}$ takes only a finite number, say $m_{k}$ of integer values on it. For $N=1,2,\ldots$, let $M_{N}$ be the product of $m_{k}$ over $k=1,2,\ldots,N$. Choose positive weights $w_{k} > 0$ such that
	\begin{enumerate}
		\item[(a)] The series $V(q) \coloneq \sum_{k=1}^{\infty}w_{k}g_{k}(q)$ converges uniformly in $q \in \mathcal{K}$;
		\item[(b)] $M_{N} R_{N} \to 0$ as $N \to \infty$, where $R_{N} = \sum_{k=N+1}^{\infty}w_{k} \|g_{k}\|_{\infty}$.
	\end{enumerate}
	It is evident that item 1) follows from (a). Indeed, $V$ is continuous by the uniform convergence. Since each $g_{k}$ is nondecreasing, $V$ is also nondecreasing. Moreover, if $p \notin \operatorname{Eq}(\mathcal{K},\varphi)$ and $t > 0$, then $g_{k}(\varphi^{t}(p)) > g_{k}(p)$ for at least one $k$. Since $w_{k} > 0$, we have that $V(\varphi^{t}(p)) > V(p)$ as well.
	
	Now we deal with item 2). If the set $\operatorname{Eq}(\mathcal{K},\varphi)$ is empty\footnote{In fact, this case cannot be realized.}, then there is nothing to prove, so we assume that it is not empty. Note that the compact $V(\operatorname{Eq}(\mathcal{K},\varphi))$ can be covered by the intervals $v + [-R_{N},R_{N}]$, where $v$ belongs to the set of possible values of the finite sum $\sum_{k=1}^{N}w_{k}g_{k}(q)$ with $q \in \operatorname{Eq}(\mathcal{K},\varphi)$. Since there are at most $M_{N}$ such values, item (b) implies item 2).
	
	Clearly, the weights 
	\begin{equation}
		w_{k} = \frac{2^{-k}}{M_{k}(1+\|g_{k}\|_{\infty})}
	\end{equation}
	satisfy both assumptions.
\end{proof}

Next, we establish that a global Lyapunov function precludes the existence of nontrivial chain recurrence in the system.
\begin{lemma}
	\label{LEM: GlobalLyapunovFunctionPreventsNontrivialRecurrence}
	Suppose that $(\mathcal{K},\varphi)$ admits a global Lyapunov function $V$ with properties as in items 1) and 2) of Lemma \ref{LEM: ZeroDimFlowGlobalFunction}. Let $q \notin \operatorname{Eq}(\mathcal{K},\varphi)$ and $T > 0$ be given. Then there exists $\varepsilon > 0$ such that there are no $(\varepsilon,T)$-chains from $q$ to $q$.
\end{lemma}
\begin{proof}
	By item 1), we have $\varkappa_{0} \coloneq V(\varphi^{T}(q)) - V(q) > 0$. By item 2), the set $\mathcal{E} \coloneq V(\operatorname{Eq}(\mathcal{K},\varphi)$ has zero measure. Thus, there exists $v \in ( V(q), V(q) + \varkappa_{0}/2]$ such that $v \notin \mathcal{E}$. Furthermore, since $\mathcal{E}$ is compact, there exists $\delta > 0$ such that $[v-\delta, v+\delta] \cap \mathcal{E} = \emptyset$. Set
	\begin{equation}
		\label{EQ: GlobalLyapFunctionChainRecurrenceFset}
		\begin{split}
		\mathcal{F} &\coloneq \{ p \in \mathcal{K} \colon | V(p) - v | \leq \delta \},\\
		\varkappa &\coloneq \min_{p \in \mathcal{F}} ( V(\varphi^{T}(p)) - V(p) ) > 0.
		\end{split}
	\end{equation}
	From the above, $\mathcal{F}$ is compact and does not intersect $\mathcal{E}$. Thus, $\varkappa$ is indeed positive, thanks to the properties of $V$.
	
	Let $\varepsilon > 0$ be such that $p_{1}, p_{2} \in \mathcal{K}$ and $d_{\mathcal{K}}(p_{1},p_{2}) < \varepsilon$ implies that
	\begin{equation}
		\label{EQ: GlobalLyapFunctionChainRecurrenceUCV}
		|V(p_{1})-V(p_{2})| \leq \min\{ \delta, \varkappa, \varkappa_{0}/2 \} \eqcolon \theta.
	\end{equation}
	
	Suppose that there exists an $(\varepsilon,T)$-chain $q_{0},q_{1},\ldots, q_{n}$ with $q_{0}=q_{n} = q$ and jump times $s_{1},\ldots, s_{n} \in [T,+\infty)$. First, we prove that if $j \in \{1,\ldots,n-1\}$ and $V(q_{j}) \geq v$, then $V(q_{j+1}) \geq v$. Suppose the opposite: $V(q_{j+1}) < v$. Using the nondecreasing of $V$, the inequality $d_{\mathcal{K}}(\varphi^{s_{j+1}}(q_{j}), q_{j+1}) < \varepsilon$ and \eqref{EQ: GlobalLyapFunctionChainRecurrenceUCV}, we obtain that
	\begin{equation}
		V(q_{j}) \leq V(\varphi^{s_{j+1}}(q_{j})) \leq V(q_{j+1}) + \theta < v + \delta.
	\end{equation}
	In particular, $q_{j}$ belongs to $\mathcal{F}$ from \eqref{EQ: GlobalLyapFunctionChainRecurrenceFset}. This gives that
	\begin{equation}
		V(\varphi^{s_{j+1}}(q_{j})) \geq V(\varphi^{T}(q_{j})) \geq V(q_{j}) + \varkappa \geq v + \varkappa.
	\end{equation}
	Combining this with the previous inequality yields
	\begin{equation}
		V(q_{j+1}) \geq v + \varkappa - \theta \geq v
	\end{equation}
	that leads to a contradiction.
	
	Using similar estimates, we obtain that
	\begin{equation}
		V(q_{1}) \geq V(\varphi^{s_{1}}(q_{0})) - \theta \geq V(\varphi^{T}(q_{0})) - \varkappa_{0}/2 \geq V(q_{0}) + \varkappa_{0}/2 \geq  v,
	\end{equation}
	so the previous reasoning can be applied by induction starting from $j=1$. It yields $V(q_{n}) = V(q) \geq v > V(q)$ that leads to a contraction.
\end{proof}

\begin{corollary}
	\label{COR: FlowsOneDimAttractors}
	Let $\varphi$ be a continuous flow on a compact metric space $\mathcal{K}$. Suppose that $\dim_{\operatorname{T}}\mathcal{K} \leq 1$ and $\check{H}^{1}(\mathcal{K};\mathbb{Z}) = 0$. Then any chain recurrent point of $\varphi$ is an equilibrium. In particular, the set $\operatorname{Eq}(\mathcal{K},\varphi)$ is not empty. Moreover, for any finite sequence of nonequilibrium points $q_{1}, \ldots, q_{n}$ such that $\omega(q_{j}) \cap \alpha(q_{j+1}) \not= \emptyset$ for any $j \in \{1,\ldots,n-1\}$, we necessarily have that $\omega(q_{n}) \cap \alpha(q_{1}) = \emptyset$. 
\end{corollary}
\begin{proof}
	The statement is obvious when $\dim_{\operatorname{T}}\mathcal{K} = 0$, since $\mathcal{K} = \operatorname{Eq}(\mathcal{K},\varphi)$. So, we assume that $\dim_{\operatorname{T}}\mathcal{K} = 1$. In this case, Lemma \ref{LEM: ZeroDimFlowGlobalFunction} ensures the existence of a global Lyapunov function, which guarantees the absence of nontrivial chain recurrence by virtue Lemma \ref{LEM: GlobalLyapunovFunctionPreventsNontrivialRecurrence}. By compactness, general arguments show that there are always exist chain recurrent points, so $\operatorname{Eq}(\mathcal{K},\varphi)$ is not empty. For the last conclusion, note that if $\omega(q_{n}) \cap \alpha(q_{1}) \not= \emptyset$, then any of the points $q_{1},\ldots,q_{n}$ is chain recurrent that contradicts the main statement.
\end{proof}

\begin{proof}[Proof of Theorem \ref{TH: BendixsonGlobalStabilityOfSemiflows}]
	By Lemma \ref{LEM: DimAndCohomInverseLimit}, we may apply Corollary \ref{COR: FlowsOneDimAttractors} to the natural extension $(\hat{\mathcal{A}},\hat{\varphi})$ of $(\mathcal{A},\varphi)$. Then Lemma \ref{LEM: ChainRecurrenceLiftInverseLimit} guarantees that the conclusion for $(\hat{\mathcal{A}},\hat{\varphi})$ can be transferred to $(\mathcal{A},\varphi)$. Finally, note that the commentary in item 2) follows from Theorem \ref{TH: BendixsonSemiflow} with $\mathcal{A} = \Gamma$.
\end{proof}
\section{Isolated fixed points of analytic mappings}
\label{SEC: IsolatedFixedPointsAnMaps}
Let $\mathbb{E}$ be a real Banach space. For any finite-dimensional subspace $\mathbb{F}$ of $\mathbb{E}$, we endow $\mathbb{F}$ with a Haar measure $\mu_{\mathbb{F}}$ normalized such that the $\mu_{\mathbb{F}}$-measure of the unit ball in $\mathbb{F}$ is the same as the Lebesgue measure of the usual unit ball in $\mathbb{R}^{k}$, where $k = \dim \mathbb{F}$. Then for any bounded linear operator $L$ in $\mathbb{E}$ and such $\mathbb{F}$, one defines the \textit{determinant} of $L$ on $\mathbb{F}$ by the formula:
\begin{equation}
	\det(L | \mathbb{F}) \coloneq \begin{cases}
		\frac{\mu_{L\mathbb{F}}(L\mathcal{B}) }{\mu_{\mathbb{F}}(\mathcal{B})}, &\text{if} \ L \ \text{is injective on} \ \mathbb{F},\\
		0, \ &\text{elsewhere},
	\end{cases}
\end{equation}
where $\mathcal{B}$ is any Borel subset of $\mathbb{F}$ with positive $\mu_{\mathbb{F}}$-measure. It can be verified that the definition does not depend on the choice of $\mathcal{B}$.

From this, for any $m=1,2,\ldots$, we define the function of singular values of $L$ by the formula:
\begin{equation}
	\omega_{m}(L) \coloneq \sup_{\dim\mathbb{F}=m} \det(L|\mathbb{F})
\end{equation}
where the infimum is taken over all closed subspaces $\mathbb{F} \subset \mathbb{E}$ with the prescribed dimension. It can be verified that $\omega_{m}(L_{2}L_{1}) \leq \omega_{m}(L_{1})\omega_{m}(L_{2})$ for any appropriately defined operators $L_{1}$ and $L_{2}$. In particular, the limit in \eqref{EQ: BanachSpaceAreaContraction} exists by Fekete's lemma.

Recall that the \textit{Kuratowski measure of noncompactness} $\alpha_{\operatorname{K}}(\mathcal{B})$ of a bounded set $\mathcal{B}$ in $\mathbb{E}$ is the infimum over $r > 0$ such that $\mathcal{B}$ can be covered by finitely many balls of radius $r$. For a bounded linear operator $L$, we set $\|L\|_{\alpha} \coloneq \alpha_{\operatorname{K}}(L\mathcal{B}_{1})$, where $\mathcal{B}_{1}$ is the unit ball in $\mathbb{E}$. It can be shown that $\|L_{2}L_{1}\|_{\alpha} \leq \|L_{1}\|_{\alpha} \|L_{2}\|_{\alpha}$ for appropriately defined operators $L_{1}$ and $L_{2}$. From this, the following limit exists:
\begin{equation}
	\rho_{\operatorname{ess}}(L) \coloneq \lim_{k \to \infty} \sqrt[k]{\|L^{k}\|_{\alpha}},
\end{equation}
which is called the \textit{essential spectral radius} of $L$.

Below, $d_{q}T$ denotes the derivative of the mapping $T$ at the point $q$. The following theorem generalizes \cite[Theorem 6]{Smith1986HD}, which, in turn, was inspired by \cite{NakajimaSeifert1983}. Similarly, to these works, the theorem can be applied to dissipative periodic differential equations.
\begin{theorem}
	\label{TH: FixedPointsOfAnaliticMapsIsolated}
	Let $\mathcal{Q}$ be a subset of a Banach space $\mathbb{E}$ with nonempty interior, and let $T$ be a continuous self-map of $\mathcal{Q}$. Suppose that
	\begin{enumerate}
		\item[1)] there exists an invariant compact $\mathcal{A}$ in the interior of $\mathcal{Q}$ that attracts points of $\mathcal{Q}$, and also attracts a superset $\mathcal{B} \supset \mathcal{A}$ with $\check{H}^{1}(\mathcal{B};\mathbb{Z}) = 0$;
		\item[2)] $T$ is analytic in a neighborhood of $\mathcal{A}$;
		\item[3)] the derivative cocycle $dT$ of $T$ over $\mathcal{A}$ uniformly contracts areas:
		\begin{equation}
			\label{EQ: BanachSpaceAreaContraction}
			\omega_{2}(dT) \coloneq \lim_{k \to \infty} \sup_{q \in \mathcal{A}} \sqrt[k]{\omega_{2}(d_{q}T^{k})} < 1.
		\end{equation}
	\end{enumerate}
	Then the fixed points of $T$ are isolated from each other. In particular, their number is finite.
\end{theorem}
\begin{proof}
	Let $\operatorname{Eq}(T)$ denote the set of fixed points of $T$. By pointwise attraction, we have that $\operatorname{Eq}(T) \subset \mathcal{A}$. Suppose that $q \in \operatorname{Eq}(T)$ is not isolated and set $L \coloneq d_{q}T$. By item 3), the essential spectral radius $\rho_{\operatorname{ess}}(L)$ of $L$ is strictly less than $1$. Indeed, if $\lambda_{1}(L) \geq \lambda_{2}(L) \geq $ are the Lyapunov exponents of $L$ counted with multiplicities, then
	\begin{equation}
		\ln \rho_{\operatorname{ess}}(L) = \inf_{k \geq 1} \lambda_{k}(L),
	\end{equation}
	see \cite[Proposition A.3.1]{Thieullen1992}.
	In particular, $\rho_{\operatorname{ess}}(L) \geq 1$ would imply that $\lambda_{k}(L) \geq 0$ for any $k$. However, $\omega_{2}(dT) < 1$ ensures that $\lambda_{1}(L) + \lambda_{2}(L) < 0$, which leads to a contradiction.
	
	Next, we note that $1$ belongs to the spectrum of $L$, since otherwise the operator $\operatorname{Id}_{\mathbb{E}} - L$ would be invertible, and the fixed point $q$ would be isolated (by the implicit function theorem). Similarly, to the above, $1$ is the only spectral point of $L$ outside the unit circle, and it is a simple eigenvalue.
	
	Let $e$ be an eigenvector of $L$ such that $Le = e$, $\mathbb{E}^{0} \coloneq \operatorname{Span}(v)$ be the corresponding spectral subspace, and let $\mathbb{E}^{s}$ be the complementary spectral subspace. Thus, $\mathbb{E} = \mathbb{E}^{s} \oplus \mathbb{E}^{0}$ is an $L$-invariant decomposition of $\mathbb{E}$. Let $P$ and $Q$ be the complementary projectors onto $\mathbb{E}^{0}$ and $\mathbb{E}^{s}$, respectively.
	
	Any $v \in \mathbb{E}$, can be uniquely decomposed as $v = q + t e + w$, where $t \in \mathbb{R}$ and $w \in \mathbb{E}^{s}$. Clearly, if $T(v) = v$ for some $v$, then $Q(v-T(v)) = 0$. Consider the function
	\begin{equation}
		\Phi(t,w) \coloneq \restr{Q(v-T(v))}{v=q + t e + w}.
	\end{equation}
	By item 2), $\Phi(t,w)$ is analytic in a neighborhood of $t=0$ and $w=0$. Moreover, the derivative of $\Phi(t,w)$ with respect to $w$ at $w=0$ is given by $\restr{(\operatorname{Id}_{\mathbb{E}} - L)}{\mathbb{E}^{s}}$, that is an invertible operator. Consequently, by the analytic implicit function theorem, there exist $\varepsilon>0$ and $\nu \colon (-\varepsilon,\varepsilon) \to \mathbb{E}^{s}$ such that $\Phi(t,w) = 0$ in a neighborhood of $t=0$ and $w=0$ if and only if $w=\nu(t)$.
	
	Let $\gamma(t) \coloneq q + t e + \nu(t)$. Clearly, $\gamma$ is an analytic function, and it is an embedding since $\gamma'(0)$ is nonzero. Any fixed point $p$ near $q$ satisfy $Q(p-T(p))=0$ and hence belongs to the image of $\gamma$. Moreover, $\gamma(t)$ is a fixed point of $T$ if and only if $g(t) \coloneq P(\gamma(t) - T(\gamma(t))) = 0$. Note that $g$ is an analytic function of one variable. Since $q$ is not isolated, $t=0$ is a limit point for zeroes of $g$. By the identity theorem for analytic functions, $g \equiv 0$. Thus, in a neighborhood of $q$, the set $\operatorname{Eq}(T)$ is parameterized by $\gamma$.
	
	Let $\mathcal{E}$ be the subset of $\operatorname{Eq}(T)$ consisting of all nonisolated points. Clearly, $\mathcal{E}$ is closed and nonempty (by our assumption). From the previous, $\mathcal{E}$ is a compact metrizable one-dimensional topological manifold without boundary. By a theorem of Milnor, $\mathcal{E}$ is a finite union of topological circles. Let $\Gamma$ be any of such circles.
	
	By item 3) and \cite[Theorem 2.5.3]{Thieullen1992}, the Hausdorff dimension of $\mathcal{A}$ is strictly less than $2$. In particular, $\dim_{\operatorname{T}}\mathcal{A} \leq 1$. This, along with item 1), justifies the applicability of Theorem \ref{TH: BendixsonMapping}, which guarantees that $\check{H}^{1}(\Gamma;\mathbb{Z}) = 0$. Thus, we arrive at a contradiction.
\end{proof}
In the context of Theorem \ref{TH: FixedPointsOfAnaliticMapsIsolated}, the set of periodic points of $T$ up to period $n$ is finite for any $n$, since the theorem can be applied to the iterates $T^{n}$ of $T$. However, the set of all periodic points may be infinite, as for the H\'{e}non map. Of course, if $T$ is the time-$1$ mapping of a continuous semiflow satisfying an analog of item 1), then Theorem \ref{TH: BendixsonGlobalStabilityOfSemiflows} is applicable.

The cohomological condition in item 1) cannot be relaxed as well. Indeed, consider the self-map of $\mathbb{R}^{2} \setminus \{(0,0)\}$ defined in polar coordinates $(r,\theta)$ by $T(r,\theta) = (\sqrt{r}, \theta)$. Here, $\mathcal{A}$ is the unit circle, which entirely consists of fixed points and attracts any closed annular region containing $\mathcal{A}$. Clearly, $\omega_{2}(d_{q}T) = 1/2$ for any $q \in \mathcal{A}$, and, consequently, item 3) is satisfied. On the other hand, $\mathcal{A}$ does not attract any superset $\mathcal{B}$ with trivial first \v{C}ech cohomology group.
\section{A refinement of Cartwright's theorems on almost periodic flows}
\label{SEC: AlmostPeriodicBendixsonRefinement}
In this section, we first aim to revisit the results of Cartwright \cite{Cartwright1967, Cartwright1969} concerning almost periodic solutions to ODEs, and present them from a more geometric perspective, using the theory of topological groups. In the case of \cite{Cartwright1967}, such a view was given by Zinchenko in \cite{Zinchenko1988}, though still in the case of ODEs. With more details, it was exposed by the present author in \cite[Chapter 4]{KuzReit2020} in Banach spaces. In \cite{MalletParet1976}, Mallet-Paret noted that the results of \cite{Cartwright1967} can be directly adapted for Hilbert spaces. As for the work \cite{Cartwright1969}, its main result is incomplete in the sense that much more can be said in the considered situation. 
\begin{remark}
	Apparently, because of their quantitative nature, these results have not received due attention in the fundamental monographs devoted to minimal flow extensions and nonautonomous dynamics (e.g., \cite{Cheban2020, ShenYi1998, LevitanZhikov1982}), where only qualitative properties are considered. It should be noted, however, that even such qualitative results acquire a much broader context of application and become linked to phenomena of a metric nature when synthesized with quantitative methods \cite{AnikushinAADyn2021, Anikushin2019Liouv, AnikushinReitRom2019, KuzReit2020}.
\end{remark}

In any case, we present self-contained proofs of the relevant generalizations of the main results from \cite{Cartwright1967, Cartwright1969} within the context of abstract minimal almost periodic flows on metric spaces (see Theorems \ref{TH: CartwrightTheoremFlows} and \ref{TH: CartwrightTheoremSkewProductFlows}), and subsequently refine them in Theorems \ref{TH: CartwrightRefinementFlows} and \ref{TH: CartwrightRefinementSkewProduct} by employing the Bendixson criterion.

To proceed further, we will need some facts from the theory of abelian topological groups \cite{Pontryagin1973}, which we will now present. Let $G$ be such a group, with its neutral element denoted by $0_{G}$. Then the \textit{character group} $\hat{G}$ of $G$ is the group of all continuous homomorphisms, called \textit{characters}, $\chi \colon G \to \mathbb{T}^{1}$ into the multiplicative circle group $\mathbb{T}^{1}$. As a subset of continuous mappings, $\hat{G}$ is endowed with the compact-open topology, which makes it an abelian topological group under the pointwise multiplication of characters. Any (continuous\footnote{In what follows, a group means an abelian topological group, any homomorphism of such groups is always assumed to be continuous, and any subgroup is assumed to be closed.}) homomorphism $i \colon G \to H$ of such groups $G$ and $H$ induces a homomorphism $i^{*} \colon \hat{H} \to \hat{G}$ of the character groups: $\chi \mapsto i^{*}\chi$, where $(i^{*}\chi)(g) = \chi(i(g))$ for $g \in G$ and $\chi \in \hat{H}$. The correspondence $G \mapsto \hat{G}$ is a contravariant functor from the category of abelian topological groups into itself. This means that $(i \circ j)^{*} = j^{*} \circ i^{*}$ for any homomorphisms $i$ and $j$ of appropriate groups, and if $i$ is the identity mapping, then so is $i^{*}$. For what follows, we need only certain results for compact\footnote{Recall that a compact group means a compact Hausdorff topological group.} or discrete metrizable groups.
\begin{description}[before=\let\makelabel\descriptionlabel]
	\item[\textbf{(TG1)}\refstepcounter{desccount}\label{DESC: TG1}] We have the following relations for compact and discrete groups:
	\begin{enumerate}
		\item[1)] If $G$ is compact, then $\hat{G}$ is discrete. If $G$ is discrete, then $\hat{G}$ is compact (\cite[Theorem 36]{Pontryagin1973});
		\item[2)] A compact group $G$ is connected if and only if $\hat{G}$ is torsion-free (\cite[Theorem 46]{Pontryagin1973});
		\item[3)] A compact group $G$ is metrizable if and only if $\hat{G}$ is countable (\cite[Theorem 43]{Pontryagin1973} and item 4) below);
		\item[4)] If $G$ is compact, then characters separate points: if $g \in G$ is such that $\chi(g) = 1$ for any $\chi \in \hat{G}$, then $g=0_{G}$ (\cite[Section 34.C)]{Pontryagin1973}).
	\end{enumerate}
\end{description}
The following key fact is known as the Pontryagin duality:
\begin{description}[before=\let\makelabel\descriptionlabel]
	\item[\textbf{(TG2)}\refstepcounter{desccount}\label{DESC: TG2}] For any compact or discrete group $G$, there is a functorial isomorphism between $G$ and the character group $\hat{\hat{G}}$ of $\hat{G}$ (\cite[Theorem 39]{Pontryagin1973}).
\end{description}
Let $H$ be a (closed) subgroup of $G$, and let $H'$ be a subgroup of $\hat{G}$. Define $\operatorname{Ann}(\hat{G},H)$ as the set of all $\chi \in \hat{G}$ such that $\chi(h) = 1$ for any $h \in H$. Similarly, let $\operatorname{Ann}(G, H')$ be the set of all $g \in G$ such that $\chi(g) = 1$ for any $\chi \in H'$. Clearly, both sets are closed subgroups.
\begin{description}[before=\let\makelabel\descriptionlabel]
	\item[\textbf{(TG3)}\refstepcounter{desccount}\label{DESC: TG3}] For any compact or discrete group $G$ and subgroups $H \subset G$ and $H' \subset \hat{G}$, there are functorial isomorphisms:
	\begin{enumerate}
		\item[1)] $\hat{H} = \hat{G}/\operatorname{Ann}(\hat{G},H)$ (\cite[Theorem 41]{Pontryagin1973});
		\item[2)] $\widehat{G/H} = \operatorname{Ann}(\hat{G}, H)$ (\cite[Theorem 37]{Pontryagin1973});
		\item[3)] $\operatorname{Ann}(G, \operatorname{Ann}(\hat{G},H)) = H$ (\cite[Theorem 40]{Pontryagin1973} or apply \nameref{DESC: TG2} with the previous items);
		\item[4)] $\operatorname{Ann}(\hat{G}, \operatorname{Ann}(G,H')) = H'$ (apply \nameref{DESC: TG2} with the previous item).
	\end{enumerate}
\end{description}

Now we turn to almost periodic functions and flows. Let $\mathcal{X}$ be a complete metric space with distance function $d_{\mathcal{X}}$, and let $\gamma \colon \mathbb{R} \to \mathcal{X}$ be a continuous $\mathcal{X}$-valued function. Given $\varepsilon > 0$, any $\tau \in \mathbb{R}$ satisfying $d_{\mathcal{X}}(\gamma(t),\gamma(t+\tau)) \leq \varepsilon$ for any $t \in \mathbb{R}$ is called an $\varepsilon$-\textit{almost period} of $\gamma$. If for any $\varepsilon > 0$ there exists $L > 0$ such that any segment $[a,a+L]$ with $a \in \mathbb{R}$ contains an $\varepsilon$-almost period of $\gamma$, then $\gamma$ is called \textit{almost periodic} (in the sense of Bohr). It is well known that any almost periodic function $\gamma$ is uniformly continuous and the closure of $\gamma(\mathbb{R})$ is compact. From this, it is clear that the property of being an almost periodic function depends only on the topological equivalence class of the distance function on $\mathcal{X}$.

By Bochner's criterion (see \cite[Section 1.2]{LevitanZhikov1982}), a bounded continuous function $\gamma$ is almost periodic if and only if the set of translates $\{\gamma(s + \cdot)\}_{s \in \mathbb{R}}$ is a precomact subset in the metric space $C_{b}(\mathbb{R}; \mathcal{X})$ of bounded continuous functions, endowed with the uniform metric. Let $\mathcal{H}(\gamma)$ be the closure of the set of translates, which is called the \textit{hull} of $\gamma$. If $\sigma$ is the translation flow on $C_{b}(\mathbb{R}; \mathcal{X})$: $\sigma^{t}(\eta) = \eta(\cdot + t)$ for $t \in \mathbb{R}$ and $\eta \in C_{b}(\mathbb{R}; \mathcal{X})$, then $(\mathcal{H}(\gamma), \sigma)$ is a minimal almost periodic flow\footnote{That is, any complete trajectory is almost periodic and dense.}. In particular, $\mathcal{H}(\gamma)$ carries a compact connected abelian group structure, for which $\gamma$ is the neutral element: if $\eta_{1},\eta_{2} \in \mathcal{H}(\gamma)$ with $\sigma^{t_{k}}(\gamma) \to \eta_{1}$ and $\sigma^{s_{k}}(\gamma) = \eta_{2}$ as $k \to \infty$, then
\begin{equation}
	\label{EQ: HullAbelianGroupStructure}
	\eta_{1} \overset{\gamma}{+} \eta_{2} \coloneq \lim_{k \to \infty} \sigma^{t_{k} + s_{k}}(\gamma).
\end{equation}
It can be shown that the definition is correct and the flow $\sigma$ agrees with the group structure: $\sigma^{t}(\eta) = \eta \overset{\gamma}{+} \sigma^{t}(\gamma)$ for any $\eta \in \mathcal{H}(\gamma)$. By \nameref{DESC: TG1}, the character group $\hat{\mathcal{H}}(\gamma)$ is a discrete countable torsion-free abelian group.
\begin{remark}
	\label{REM: ApModulesMetricSpaces}
	If $\mathcal{X}$ is a Banach space, then there is the Fourier theory producing exponents $\lambda_{1},\lambda_{2},\ldots$ of $\gamma$, see \cite[Chapter 2]{LevitanZhikov1982}. It can be shown that the smallest $\mathbb{Z}$-module in $\mathbb{R}$ containing the exponents is algebraically isomorphic to $\hat{\mathcal{H}}(\gamma)$, see \cite{Zinchenko1988} or \cite[Theorem 4.7]{KuzReit2020}. In fact, this module is $\mathfrak{M}_{\mathbb{Z}}(\hat{\mathcal{H}}(\gamma), \sigma)$ defined in \eqref{EQ: ZModuleAPFlowDef}.
\end{remark}
\begin{remark}
	\label{REM: NonintersectionOfAp}
	Note that if $\gamma_{1}$ and $\gamma_{2}$ are almost periodic functions for which there exists $t_{0} \in \mathbb{R}$ such that $\gamma_{1}(t) =\gamma_{2}(t)$ for any $t \geq t_{0}$, then $\gamma_{1}(t) = \gamma_{2}(t)$ for all $t \in \mathbb{R}$. Indeed, by Bochner's criterion, $t \mapsto (\gamma_{1}(t), \gamma_{2}(t)) \in \mathcal{X} \times \mathcal{X}$ is almost periodic. Thus, for any $\varepsilon > 0$ there exist arbitrarily large common $\varepsilon$-almost periods $\tau$. In particular,
	\begin{equation}
		\begin{split}
			d_{\mathcal{X}}(\gamma_{1}(t), \gamma_{2}(t)) \leq \\ d_{\mathcal{X}}(\gamma_{1}(t), \gamma_{1}(t+\tau)) + d_{\mathcal{X}}(\gamma_{1}(t+\tau), \gamma_{2}(t+\tau)) + d_{\mathcal{X}}(\gamma_{2}(t), \gamma_{2}(t+\tau)) \leq 2\varepsilon.
		\end{split}
	\end{equation}
	for any $t$ such that $t + \tau \geq t_{0}$. Since $\varepsilon>0$ can be arbitrarily small, and $\tau$ can be arbitrarily large, we must have the equality $\gamma_{1}(t) = \gamma_{2}(t)$ for all $t \in \mathbb{R}$. In particular, almost periodic trajectories of a semiflow cannot intersect.
\end{remark}
Next, we consider an abstract context. Let $(G,+,0_{G})$ be a metrizable compact connected abelian group with group operation $+$ and neutral element $0_{G}$. Suppose that $\varphi$ is a continuous flow on $G$ that agrees with the group structure: $\varphi^{t}(g) = g + \varphi^{t}(0_{G})$ for any $t \in \mathbb{R}$ and $g \in G$. We call $\varphi$ or, in more detail, the tuple $((G,+,0_{G}), \varphi)$ a $G$-\textit{consistent flow}. When $\varphi$ is initially defined on a metrizable compact space $\mathcal{X}$, and $G$ is a group structure on $\mathcal{X}$ that makes $\varphi$ a $G$-consistent flow, we also say that $G$ is $\varphi$-\textit{consistent}. 
\begin{remark}
	By $G$-consistency, $\varphi$ is two-sided uniformly Lyapunov stable. By \cite[Property 8, Section 1.3]{LevitanZhikov1982}, it is almost periodic (not necessarily minimal).
\end{remark}
Consider the homomorphism $i_{\varphi} \colon \mathbb{R} \to G$ given by $i_{\varphi}(t) = \varphi^{t}(0_{G})$. It induces a homomorphism on the character groups $i^{*}_{\gamma} \colon \hat{G} \to \hat{\mathbb{R}}$. It is well known that any character on $\mathbb{R}$ has the form $\chi_{\lambda}(t) = e^{i \lambda t}$ for some $\lambda \in \mathbb{R}$. Define the \textit{frequency $\mathbb{Z}$-module} of $\varphi$ as follows:
\begin{equation}
	\label{EQ: ZModuleAPFlowDef}
	\mathfrak{M}_{\mathbb{Z}}(\varphi) = \mathfrak{M}_{\mathbb{Z}}(G, \varphi) \coloneq \{ \lambda \in \mathbb{R} \ | \ \chi_{\lambda} \in i^{*}_{\varphi}(\hat{G}) \}.
\end{equation}
By \nameref{DESC: TG1}, this is a countable $\mathbb{Z}$-submodule in $\mathbb{R}$. Note that $\varphi$ is minimal if and only if $i_{\varphi}$ has dense range. In this case, $i^{*}_{\varphi}$ is injective and hence induces an algebraic isomorphism between $\hat{G}$ and $\mathfrak{M}_{\mathbb{Z}}(\varphi)$.

\begin{remark}
	In what follows, when speaking of such $G$, we treat it either as a group structure or as a set of group elements, depending on which is more appropriate in the given context.
\end{remark}
\begin{remark}
	It will be also convenient to consider the \textit{frequency $\mathbb{Q}$-module} $\mathfrak{M}_{\mathbb{Q}}(\varphi) = \mathfrak{M}_{\mathbb{Z}}(\varphi) \otimes \mathbb{Q}$ of $\varphi$, i.e., the smallest $\mathbb{Q}$-module in $\mathbb{R}$ containing $\mathfrak{M}_{\mathbb{Z}}(\varphi)$. Since $\mathfrak{M}_{\mathbb{Z}}(\varphi)$ is torsion-free, we have that $\operatorname{rank}\mathfrak{M}_{\mathbb{Z}}(\varphi) = \dim \mathfrak{M}_{\mathbb{Q}}(\varphi)$.
\end{remark}

A $G$-consistent flow $\varphi$ and an $H$-consistent flow $\psi$ are said to be \textit{isomorphic} if there is a topological isomorphism of $G$ and $H$ which additionally conjugates $\varphi$ and $\psi$.
\begin{proposition}
	\label{PROP: GConsistencFlowsFrequencyModule}
	Let $\varphi$ and $\psi$ be isomorphic $G$-consistent and $H$-consistent flows, respectively. Then $\mathfrak{M}_{\mathbb{Z}}(\varphi) = \mathfrak{M}_{\mathbb{Z}}(\psi)$.
\end{proposition}
\begin{proof}
	Let $\alpha \colon G \to H$ be an isomorphism that cojugates the flows. Then $\alpha(i_{\varphi}(t)) = \alpha(\varphi^{t}(0)) = \psi^{t}(0) = i_{\psi}(t)$ for any $t \in \mathbb{R}$. Thus, $i_{\psi} = \alpha \circ i_{\varphi}$ and hence $i^{*}_{\psi} = i^{*}_{\varphi} \circ \alpha^{*}$. Since $\alpha$ is an isomorphism of topological groups, $\alpha^{*}$ is an isomorphism of the character groups. Consequently, $i^{*}_{\psi}(\hat{H}) = i^{*}_{\varphi}(\hat{G})$, as desired.
\end{proof}

Let $(\mathcal{X},\varphi)$ be a minimal almost periodic flow. Similarly to \eqref{EQ: HullAbelianGroupStructure}, any $x \in \mathcal{X}$ determines on $\mathcal{X}$ a $\varphi$-consistent group structure $G$ with group operation $\overset{x}{+}$ and neutral element $x$. For distinct points $x,y \in \mathcal{X}$, the translation $z \mapsto z \overset{x}{-} y \overset{x}{+}x$ provides an isomorphism between $((\mathcal{X}, \overset{y}{+}, y), \varphi)$ and $((\mathcal{X}, \overset{x}{+}, x), \varphi)$. Consequently, $\mathfrak{M}_{\mathbb{Z}}(\varphi)$ does not depend on the choice of $x$, thanks to Proposition \ref{PROP: GConsistencFlowsFrequencyModule}. Note that any $\varphi$-consistent group structure can be obtained in this way.
\begin{remark}
	If $\gamma$ is an almost periodic trajectory of $(\mathcal{X},\varphi)$, then it is easily seen that the correspondence $\mathcal{H}(\gamma) \ni \eta \mapsto \eta(0) \in \mathcal{X}$ determines an isomorphism between $(\mathcal{H}(\gamma), \overset{\gamma}{+}, \sigma)$ and $(\mathcal{X}, \overset{x}{+}, \varphi)$, where $x=\gamma(0)$. 
\end{remark}

More generally, we have the following.
\begin{theorem}
	\label{TH: ConjugationOfMinAPFlows}
	Let $(\mathcal{X}, \varphi)$ and $(\mathcal{Y}, \psi)$ be minimal almost periodic flows. Then the flows are topologically conjugate if and only if $\mathfrak{M}_{\mathbb{Z}}(\varphi) = \mathfrak{M}_{\mathbb{Z}}(\psi)$.
\end{theorem}
\begin{proof}
	The ``only if'' part follows just from the above considerations. Specifically, we take any $x \in \mathcal{X}$ and $y \in \mathcal{Y}$ related by the conjugacy and consider the corresponding group structures. Then the conjugacy becomes an isomorphism of the consistent flows. Thus, Proposition \ref{PROP: GConsistencFlowsFrequencyModule} guarantees the coincidence of frequency modules.
	
	The ''if'' part is based on the Pontryagin duality. Fix any $x \in \mathcal{X}$ and $y \in \mathcal{Y}$. Let $i_{\varphi}(t) = \varphi^{t}(x)$ and $i_{\psi}(t) = \psi^{t}(y)$ be the homomorphisms of $\mathbb{R}$ into the corresponding groups $G \leftrightarrow (\mathcal{X},\overset{x}{+}, x)$ and $H \leftrightarrow (\mathcal{Y},\overset{y}{+}, y)$. By minimality, the restrictions $i^{*}_{\varphi} \colon \hat{G} \to \mathfrak{M}_{\mathbb{Z}}(\varphi)$ and $i^{*}_{\psi} \colon \hat{H} \to \mathfrak{M}_{\mathbb{Z}}(\psi)$ are isomorphisms. Since the frequency modules coincide, there is an isomorphism $j^{*} \coloneq (i^{*}_{\varphi})^{-1} \circ i^{*}_{\psi} \colon \hat{H} \to \hat{G}$. By \nameref{DESC: TG2}, it induces an isomorphism $j \colon G \to H$ of topological groups. Clearly, $j (\varphi^{t}(x)) =  \psi^{t}(y)$ for any $t \in \mathbb{R}$. Finally, by consistency, we have that
	\begin{equation}
		j( \varphi^{t}(z) ) = j( z \overset{x}{+} \varphi^{t}(x) ) = j(z) \overset{y}{+} \psi^{t}(y) = \psi^{t}(j(z))
	\end{equation}
	for any $z \in \mathcal{X}$ and $t \in \mathbb{R}$. Thus, $\varphi$ and $\psi$ are topologically conjugate.
\end{proof}
\begin{remark}
	Similar considerations show that $(\mathcal{X}, \varphi)$ is an extension of $(\mathcal{Y}, \psi)$ if and only if $\mathfrak{M}_{\mathbb{Z}}(\varphi) \supset \mathfrak{M}_{\mathbb{Z}}(\psi)$.
\end{remark}

A certain general version of the main theorem of Cartwright from \cite{Cartwright1967} can be stated as follows.
\begin{theorem}
	\label{TH: CartwrightTheoremFlows}
	Let $\varphi$ be a continuous semiflow on a complete metric space $\mathcal{X}$, and suppose that $\gamma$ is an almost periodic trajectory of $\varphi$. Let $\Gamma$ be the closure of $\gamma$. Then 
	\begin{enumerate}
		\item[1)] $(\Gamma,\varphi)$ is a minimal almost periodic flow\footnote{More rigorously, $\varphi$ can be uniquely extended to such a flow.};
		\item[2)] $\dim_{\operatorname{T}}\Gamma = \operatorname{rank}\mathfrak{M}_{\mathbb{Z}}(\Gamma,\varphi)$.
	\end{enumerate}
	Moreover, $(\Gamma,\varphi)$ is $T$-periodic ($T > 0$) if and only if $\mathfrak{M}_{\mathbb{Z}}(\Gamma,\varphi) \subset \frac{2\pi}{T}\mathbb{Z}$, i.e., $\mathfrak{M}_{\mathbb{Z}}(\Gamma,\varphi) = \frac{2\pi N}{T}\mathbb{Z}$ for some $N \in \mathbb{Z}$.
\end{theorem}
\begin{proof}
	Item 1) can be shown similarly to \cite{Cartwright1967} or \cite[Theorems 4.4 and 4.5]{KuzReit2020}, taking into account Remark \ref{REM: NonintersectionOfAp}. In view of the group structure on $\Gamma$ and the definition of $\mathfrak{M}_{\mathbb{Z}}(\Gamma,\varphi)$, item 2) becomes a general fact from the theory of topological groups, see \cite[Theorem 47]{Pontryagin1973} or \cite[Theorem 4.2]{KuzReit2020}.
	
	For the last part, we endow $\Gamma$ with any $\varphi$-consistent group structure $G$ and recall that $i_{\varphi}(t) = \varphi^{t}(0_{G})$. Clearly, $\varphi$ is $T$-periodic if and only if $i_{\varphi}(T) = \varphi^{T}(0_{G}) = 0_{G}$. Since characters separate points (see \nameref{DESC: TG1}), this is equivalent to that $e^{i\lambda T} = 1$ for any $\lambda \in \mathfrak{M}_{\mathbb{Z}}(\Gamma,\varphi)$. This shows that the frequency module is a subgroup of $\frac{2\pi}{T}\mathbb{Z}$. Since any such subgroup must be cyclic, the statement is proved.
\end{proof}

When $\mathcal{X} = \mathbb{R}^{n}$ in the above theorem, we must necessarily have $\dim_{\operatorname{T}}\Gamma \leq n-1$, since compact minimal sets cannot have full dimension\footnote{This is because the boundary of an invariant set is invariant, and sets with full topological dimension necessarily have nonempty interior.}. Moreover, if $\dim_{\operatorname{T}}\Gamma = n-1$, then any $\varphi$-consistent group structure $G$ on $\Gamma$ is topologically isomorphic to the torus group, see \cite{Cartwright1967}. Using the Bendixson criterion, we can refine the former conclusion in terms of dimension estimates and cohomology.
\begin{theorem}
	\label{TH: CartwrightRefinementFlows}
	In the context of Theorem \ref{TH: CartwrightTheoremFlows}, suppose that there exist a positive integer $m$ and a compact set $\mathcal{A} \supset \Gamma$ satisfying the following:
	\begin{enumerate}
		\item[1)] $\mathcal{A}$ attracts a set $\mathcal{B} \supset \Gamma$ with $\check{H}^{m}(\mathcal{B};\mathbb{Z}) = 0$;
		\item[2)] $\dim_{\operatorname{T}}\mathcal{A} \leq m$ (if $\Gamma \not= \mathcal{A}$).
	\end{enumerate}
	Then $\dim_{\operatorname{T}}\Gamma = \operatorname{rank}\mathfrak{M}_{\mathbb{Z}}(\Gamma,\varphi) \leq m - 1$.
\end{theorem}
\begin{proof}
	By item 2), we have $\dim_{\operatorname{T}}\Gamma \leq m$. By endowing $\Gamma$ with any $\varphi$-consistent group structure, we can conclude that if $\dim_{\operatorname{T}}\Gamma = m$, then $\check{H}^{m}(\Gamma;\mathbb{Z}) \neq 0$, as we discussed at the end of Section \ref{SEC: CohomologicalBendixsonCriterion}. On the other hand, Theorem \ref{TH: BendixsonSemiflow} ensures that $\check{H}^{m}(\Gamma;\mathbb{Z}) = 0$, which leads to a contradiction. Consequently, $\dim_{\operatorname{T}}\Gamma \leq m-1$.
\end{proof}

Next, we turn to the main result of \cite{Cartwright1969}. Before stating its appropriate generalization, we introduce some preliminary notations. Let $\wp \colon (\mathcal{X},\pi) \to (\mathcal{Q},\vartheta)$ be an extension, where both $(\mathcal{X},\pi)$ and $(\mathcal{Q},\vartheta)$ are minimal almost periodic flows. Choosing $x_{0} \in \mathcal{X}$ and $q_{0}=\wp(x_{0})$, we endow $\mathcal{X}$ and $\mathcal{Q}$ with group structures such that $\wp$ becomes a homomorphism. Then for $i_{\pi}(t)=\pi^{t}(x_{0})$ and $i_{\vartheta}(t)=\vartheta^{t}(q_{0})$, it follows that $i^{*}_{\vartheta} = i^{*}_{\pi} \circ \wp^{*}$. In particular, we have the inclusion $\mathfrak{M}_{\mathbb{Z}}(\vartheta) \subset \mathfrak{M}_{\mathbb{Z}}(\pi)$. 

Proceeding in this context, let $\operatorname{Tor}(\pi, \vartheta)$ be the torsion subgroup of the frequency quotient $\mathfrak{M}_{\mathbb{Z}}(\pi)/\mathfrak{M}_{\mathbb{Z}}(\vartheta)$, i.e., the subgroup consisting of all elements with finite order. Let
\begin{equation}
	\operatorname{Sat}(\pi,\vartheta) = \mathfrak{M}_{\mathbb{Z}}(\pi) \cap \mathfrak{M}_{\mathbb{Q}}(\vartheta) = \{ \lambda \in \mathfrak{M}_{\mathbb{Z}}(\pi) \ | \ n \lambda \in \mathfrak{M}_{\mathbb{Z}}(\vartheta) \ \text{for some} \ n \in \mathbb{Z} \}.
\end{equation}  
be the $\mathbb{Z}$-saturation of $\mathfrak{M}_{\mathbb{Z}}(\vartheta)$ in $\mathfrak{M}_{\mathbb{Z}}(\pi)$. Then it is not hard to see that
\begin{equation}
	\label{EQ: CartSkewProdTorsViaSatQuotient}
	\operatorname{Tor}(\pi, \vartheta) = \operatorname{Sat}(\pi,\vartheta)/\mathfrak{M}_{\mathbb{Z}}(\vartheta).
\end{equation}
The elements of $\operatorname{Tor}(\pi, \vartheta)$ correspond to internal resonances $\lambda \in \mathfrak{M}_{\mathbb{Z}}(\pi)$, namely, such that $n \lambda \in \mathfrak{M}_{\mathbb{Z}}(\vartheta)$ for some $n \geq 2$, but $\lambda \notin \mathfrak{M}_{\mathbb{Z}}(\vartheta)$.

For the torsion-free quotient, we have the following standard relation:
\begin{equation}
	\label{EQ: CartSkewProdTorsFreeQuotViaSat}
	\left(\mathfrak{M}_{\mathbb{Z}}(\pi)/\mathfrak{M}_{\mathbb{Z}}(\vartheta)\right)/\operatorname{Tor}(\pi, \vartheta) = \mathfrak{M}_{\mathbb{Z}}(\pi)/\operatorname{Sat}(\pi,\vartheta).
\end{equation}

For a compact abelian group $K$, let $K_{0}$ denote the identity component of $K$, i.e., the connected component containing the neutral element. It is well known that $K_{0}$ is a (closed) subgroup of $K$, and the quotient $K/K_{0}$ is totally disconnected and can be identified with the set of connected components $\pi_{0}(K)$ of $K$. We will need the following well-known fact, for which we provide a proof, as similar reasoning will be used later.
\begin{proposition}
	\label{PROP: CompactGroupUnitComponent}
	For $K$ and $K_{0}$ as above, we have that $\operatorname{Ann}(\hat{K},K_{0}) = \operatorname{Tor}(\hat{K})$ and $K/K_{0} = \operatorname{Hom}(\operatorname{Tor}(\hat{K}), \mathbb{T}^{1})$.
\end{proposition}
\begin{proof}
	Set $K' \coloneq \operatorname{Ann}(K, \operatorname{Tor}(\hat{K}))$. By item 4) of \nameref{DESC: TG3}, we have that $\operatorname{Ann}(\hat{K},K') = \operatorname{Tor}(\hat{K})$. This and item 1) of \nameref{DESC: TG3} ensure that $\hat{K}' = \hat{K} / \operatorname{Tor}(\hat{K})$. By item 2) of \nameref{DESC: TG1}, $K'$ is connected. By the definition of $K_{0}$, this gives $K' \subset K_{0}$. On the other hand, since $K_{0}$ is connected, $\hat{K}_{0} = \hat{K}/\operatorname{Ann}(\hat{K}, K_{0})$ is torsion-free. In particular, the inclusion $\operatorname{Ann}(\hat{K}, K_{0}) \supset \operatorname{Tor}(\hat{K})$ is valid. Passing to the annihilators in $K$, we obtain that $K_{0} \subset K'$ and hence $K_{0} = K'$. 
	
	Using this and item 2) of \nameref{DESC: TG3} gives the following relation:
	\begin{equation}
		\widehat{K/K_{0}} = \operatorname{Tor}(\hat{K}).
	\end{equation}
	Passing to the character groups (note that $\operatorname{Tor}(\hat{K})$ is discrete) and using \nameref{DESC: TG2}, we obtain the last conclusion.
\end{proof}

In the following theorem, we present a generalization of the main result from \cite{Cartwright1969}.
\begin{theorem}
	\label{TH: CartwrightTheoremSkewProductFlows}
	Let $(\mathcal{X},\pi)$ and $(\mathcal{Q},\vartheta)$ be minimal almost periodic flows, with the former being an extension of the latter via a mapping $\wp \colon \mathcal{X} \to \mathcal{Q}$. Then there exists an almost periodic flow $\varphi$ on $\mathcal{X}$ that preserves the fibers ($\wp = \wp \circ \varphi^{s}$ for any $s \in \mathbb{R}$) and is consistent with the skew-product motion:
	\begin{equation}
		\label{EQ: CartwrightSkewProductMotionConsistency}
		\varphi^{s}(\pi^{t}(x)) = \pi^{t}(\varphi^{s}(x)) \qquad \text{for any} \, x \in \mathcal{X}, \ \text{and} \ t,s \in \mathbb{R}.
	\end{equation}
	For any minimal component of $\varphi$, there is a $\pi$-consistent group structure $G$ on $\mathcal{X}$ such that the component is of the form $(K_{0},\varphi)$, where $K$ is the fiber (subgroup) containing $0_{G}$. In this context, we have the following:
	\begin{enumerate}
		\item[1)] Other minimal components of $\varphi$ are topologically conjugate to $(K_{0},\varphi)$ by appropriate translations on $G$;
		\item[2)] $\mathfrak{M}_{\mathbb{Z}}(K_{0},\varphi) \cong \mathfrak{M}_{\mathbb{Z}}(\pi)/\operatorname{Sat}(\pi,\vartheta)$ \quad and \quad $\mathfrak{M}_{\mathbb{Q}}(K_{0},\varphi) \cong \mathfrak{M}_{\mathbb{Q}}(\pi)/\mathfrak{M}_{\mathbb{Q}}(\vartheta)$;
		\item[3)] $\pi_{0}(K) \coloneq K/K_{0} \cong \operatorname{Hom}(\operatorname{Tor}(\pi, \vartheta); \mathbb{T}^{1})$.
	\end{enumerate}
\end{theorem}
\begin{remark}
	\label{REM: VerticalFlow}
	In the context of Theorem \ref{TH: CartwrightTheoremSkewProductFlows}, when $G$ is fixed, then any $\varphi$ satisfying \eqref{EQ: CartwrightSkewProductMotionConsistency} is $G$-consistent. In particular, it is uniquely determined by the trajectory $s \mapsto \varphi^{s}(0_{G})$, thanks to the density of $t \mapsto \pi^{t}(0_{G})$ in $\mathcal{X}$. Moreover, the isomorphisms in item 2) are functorial in the category whose objects are tuples $((G,\pi), (H,\vartheta), \wp, \varphi)$, with morphisms defined in an appropriate manner. For the functoriality of the isomorphism in item 3), $\varphi$ is not needed.
\end{remark}
\begin{proof}
	Similarly to the above, we fix points $x_{0} \in \mathcal{X}$ and $q_{0} = \wp(x_{0})$ and consider the corresponding group structures $G$ and $H$ on $\mathcal{X}$ and $\mathcal{Q}$, respectively. For the associated homomorphisms $i_{\pi} \colon \mathbb{R} \to G$ and $i_{\vartheta} \colon \mathbb{R} \to H$ given by $i_{\pi}(t) = \pi^{t}(x_{0})$ and $i_{\vartheta}(t) = \vartheta^{t}(q_{0})$, we have that $\wp \circ i_{\pi} = i_{\vartheta}$. Since the ranges of $i_{\pi}$ and $i_{\vartheta}$ are dense, $\wp \colon G \to H$ becomes a surjective homomorphism of topological groups. In particular, $\chi \mapsto \wp^{*}\chi \coloneq \chi \circ \wp$ embeds $\hat{H}$ into $\hat{G}$. Clearly, $\mathcal{X}_{q_{0}} = \wp^{-1}(q_{0}) = \operatorname{Ker}\wp$, and for any $q \in \mathcal{Q}$, the fiber $\mathcal{X}_{q} = \wp^{-1}(q)$ is the coset in $G/\operatorname{Ker}\wp$. Since $G/\operatorname{Ker}\wp \cong H$, items 2) and 3) of \nameref{DESC: TG3} give that
	\begin{equation}
		\label{EQ: CartwrightKerPdescr}
		\operatorname{Ker}\wp = \operatorname{Ann}(G, \wp^{*}(\hat{H})) \coloneq \{g \in G \ | \ (\wp^{*}\chi)(g) = 1 \ \text{for any} \ \chi \in \hat{H}\}.
	\end{equation}
	
	Next, we are going to construct an almost periodic trajectory $\mathbb{R} \ni t \mapsto f(t) \in K \coloneq \operatorname{Ker}\wp$ using duality. Since $\operatorname{Ann}(\hat{G}, K) = \wp^{*}(\hat{H})$ by the same principles, using item 1) of \nameref{DESC: TG3} ensures that
	\begin{equation}
		\label{EQ: CartSPCharKformula}
		\hat{K} \cong \hat{G}/\wp^{*}(\hat{H}) 
		\cong \mathfrak{M}_{\mathbb{Z}}(\pi)/\mathfrak{M}_{\mathbb{Z}}(\vartheta),
	\end{equation}
	where in the last isomorphism the character groups are identified with the corresponding frequency modules according to $i^{*}_{\pi}$ and $i^{*}_{\vartheta}$, respectively.
	
	Next, we fix any $\mathbb{Q}$-linear complement $V$ to the subspace $\mathfrak{M}_{\mathbb{Q}}(\vartheta)$ in $\mathfrak{M}_{\mathbb{Q}}(\pi)$, so $\mathfrak{M}_{\mathbb{Q}}(\pi) = V \oplus \mathfrak{M}_{\mathbb{Q}}(\vartheta)$. Let $P$ be the projector onto $V$. Then $\theta = \restr{P}{\mathfrak{M}_{\mathbb{Z}}(\pi)} \colon \mathfrak{M}_{\mathbb{Z}}(\pi) \to V \subset \mathbb{R}$ is a homomorphism with $\operatorname{Ker}\theta = \mathfrak{M}_{\mathbb{Z}}(\pi) \cap \mathfrak{M}_{\mathbb{Q}}(\vartheta) = \operatorname{Sat}(\pi,\vartheta)$. For any $t \in \mathbb{R}$, define
	\begin{equation}
		\label{EQ: CartwrightSPDenseOrbitDef}
		f(t) \coloneq ( \lambda \mapsto e^{i \theta(\lambda) t} ) \in \operatorname{Hom}(\mathfrak{M}_{\mathbb{Z}}(\pi); \mathbb{T}^{1}) \cong \operatorname{Hom}(\hat{G}; \mathbb{T}^{1}) \cong G,
	\end{equation}
	where in the last equality we used the duality, i.e., \nameref{DESC: TG2}. Since $\hat{G}$ is discrete, the compact-open topology on $\operatorname{Hom}(\hat{G}; \mathbb{T}^{1})$ is the metrizable topology of pointwise convergence. Thus, $f \colon \mathbb{R} \to G$ is a continuous homomorphism. By construction and \eqref{EQ: CartwrightKerPdescr}, we have that $f(t) \in K = \operatorname{Ker}\wp$ for any $t \in \mathbb{R}$. Consequently, the closure $K_{f}$ of $f(\mathbb{R})$ is a subgroup of $K$.
	
	Now define a global flow $\varphi$ on $\mathcal{X} = G$ by setting $\varphi^{t}(g) = g + f(t)$ for any $g \in G$ and $t \in \mathbb{R}$, where $+ \coloneq \overset{x_{0}}{+}$. Since $f$ is a homomorphism, for any $t,s \in \mathbb{R}$ and $g \in G$, we have that
	\begin{equation}
		\varphi^{t+s}(g) = g + f(t+s) = g + f(s) + f(t) = \varphi^{t}(g+f(s)) = \varphi^{t}(\varphi^{s}(g)),
	\end{equation}
	that is, $\varphi$ is indeed a continuous flow. Clearly, any motion $t \mapsto \varphi^{t}(g)$ is almost periodic and lies in the coset of $g$ in $G/K$. Since $\pi$ is $G$-consistent, we obtain that
	\begin{equation}
		\begin{split}
			\varphi^{s}(\pi^{t}(x)) = \pi^{t}(x) + f(s) = x + f(s) + \pi^{t}(x_{0}),\\
			\pi^{t}(\varphi^{s}(x)) = \pi^{t}(x + f(s)) = x + f(s) + \pi^{t}(x_{0})
		\end{split}
	\end{equation}
	for any $x \in \mathcal{X}$ and $t,s \in \mathbb{R}$. This shows \eqref{EQ: CartwrightSkewProductMotionConsistency}.
	
	By Theorem \ref{TH: CartwrightTheoremFlows}, $(K_{f}, \varphi)$ is a minimal almost periodic flow. Since $\operatorname{Ann}(\hat{G},K_{f}) \cong \operatorname{Ker}\theta$ and $\operatorname{Ker}\theta = \operatorname{Sat}(\pi,\vartheta)$, we have that
	\begin{equation}
		\label{EQ: CartSkewProdCharKf}
		\hat{K}_{f} \cong \hat{G}/\operatorname{Ann}(\hat{G}, K_{f}) \cong \mathfrak{M}_{\mathbb{Z}}(\pi)/\operatorname{Sat}(\pi,\vartheta) \cong \theta(\mathfrak{M}_{\mathbb{Z}}(\pi)).
	\end{equation}
	In particular, $\mathfrak{M}_{\mathbb{Z}}(K_{f}, \varphi) = \theta(\mathfrak{M}_{\mathbb{Z}}(\pi)) \cong \mathfrak{M}_{\mathbb{Z}}(\pi)/\operatorname{Sat}(\pi,\vartheta)$. Moreover,
	\begin{equation}
		\label{EQ: CarthwrightSkewProductFrequenciesProof}
		\mathfrak{M}_{\mathbb{Q}}(K_{f}, \varphi) = P(\mathfrak{M}_{\mathbb{Q}}(\pi)) = V \cong \mathfrak{M}_{\mathbb{Q}}(\pi) / \mathfrak{M}_{\mathbb{Q}}(\vartheta).
	\end{equation}
	
	Let $g \in G$ belong to a minimal component of $\varphi$, say $K_{g}$. By virtue of the above, it is not hard to see that the translation $h \mapsto h - g$ provides a topological conjugacy between $(K_{g},\varphi)$ and $(K_{f}, \varphi)$. Since $\varphi$ is $G$-consistent with respect to any $\pi$-consistent group structure $G$ (see Remark \ref{REM: VerticalFlow}), the same argument can be applied when defining $G$ using any minimal component of $\varphi$. This, combined with Theorem \ref{TH: ConjugationOfMinAPFlows}, proves items 1) and 2).
	
	Now we turn to item 3). Let $K_{0}$ be the connected component of the neutral element $x_{0}$. Since $x_{0}=f(0) \in K_{f}$ and $K_{f}$ is connected, we have that $K_{f} \subset K_{0}$. By \eqref{EQ: CartwrightSPDenseOrbitDef}, \eqref{EQ: CartSkewProdCharKf}, and \eqref{EQ: CartSkewProdTorsViaSatQuotient}, we have that 
	\begin{equation}
		\operatorname{Ann}(\hat{K}, K_{f}) \cong \operatorname{Sat}(\pi,\vartheta)/\mathfrak{M}_{\mathbb{Z}}(\vartheta) \cong \operatorname{Tor}(\pi,\vartheta) \cong \operatorname{Tor}(\hat{K}).
	\end{equation}
	Using Proposition \ref{PROP: CompactGroupUnitComponent} and item 3) of \nameref{DESC: TG3}, we obtain that $K_{0} = K_{f}$ and that other connected components in the fiber $K=\operatorname{Ker}\wp=\wp^{-1}(q_{0})$ correspond to the cosets in $K / K_{f} = K / K_{0} \cong \operatorname{Hom}(\operatorname{Tor}(\pi,\vartheta); \mathbb{T}^{1})$.
\end{proof}
\begin{remark}
	In terms of general nonautonomous dynamics \cite{LevitanZhikov1982, Cheban2020, ShenYi1998}, the choice of a group structure $G$ provides an isomorphism with the enveloping Ellis semigroup $E(\mathcal{X},\pi)$: $T \mapsto T(0_{G})$. In the case of minimal almost periodic flows, this group acts freely and transitively that allows to identify the fiber containing $0_{G}$ with a subgroup, thereby expressing dynamical and topological properties in simple algebraic terms. In this sense, items 2) and 3) of Theorem \ref{TH: CartwrightTheoremSkewProductFlows} describe the structure of the Ellis semigroup of an extension. Beyond almost periodic systems, the situation is much more complicated, see Example \ref{EX: DenjoyAlmostAutomorphic} below.
\end{remark}

In the context of Theorem \ref{TH: CartwrightTheoremSkewProductFlows}, we call the flow $\varphi$ a \textit{vertical flow} of the extension $\wp \colon (\mathcal{X},\pi) \to (\mathcal{Q},\vartheta)$. Any minimal component of $\varphi$ is subject to possible restrictions arising from the embedding of fibers $\mathcal{X}_{q}$ in $\mathbb{R}^{n}$, discussed after Theorem \ref{TH: CartwrightTheoremFlows}. Similarly to Theorem \ref{TH: CartwrightRefinementFlows}, we can refine these constraints using the Bendixson criterion as follows.
\begin{theorem}
	\label{TH: CartwrightRefinementSkewProduct}
	Let $\wp \colon (\mathcal{Z},\pi) \to (\mathcal{Q},\vartheta)$ be an extension of a minimal almost periodic flow $(\mathcal{Q},\vartheta)$ to a continuous semiflow $(\mathcal{Z},\pi)$ in a complete metric space $\mathcal{Z}$. Suppose that $\gamma$ is an almost periodic trajectory of $\pi$, and let $\mathcal{X}$ be the closure of $\gamma$. Then $(\mathcal{X},\pi)$ is a minimal almost periodic flow that is an extension of $(\mathcal{Q},\vartheta)$ via the mapping $\restr{\wp}{\mathcal{X}}$, thereby satisfying the assumptions of Theorem \ref{TH: CartwrightTheoremSkewProductFlows}. Let $(\mathcal{X},\varphi)$ be a vertical flow of the corresponding extension.
	
	Proceeding in the above context, suppose that there exist a positive integer $m$, $q \in \mathcal{Q}$, and a compact set $\mathcal{A} \supset \Gamma \coloneq \wp^{-1}(q) \cap \mathcal{X}_{q}$ in the fiber $\mathcal{Z}_{q}$ over $q$ satisfying the following:
	\begin{enumerate}
		\item[1)] $\mathcal{A}$ pullback attracts over $q$ a family of sets $\{ \mathcal{B}_{s} \}_{s \leq 0}$ such that $\Gamma_{s} \subset \mathcal{B}_{s} \subset \mathcal{Z}_{\vartheta^{s}(q)}$, where $\Gamma_{s} \coloneq (\pi^{-s})^{-1}(\Gamma) \cap \mathcal{X}_{\vartheta^{s}(q)}$, and $\check{H}^{m}(\mathcal{B}_{s};\mathbb{Z}) = 0$ for any $s \leq 0$;
		\item[2)] $\dim_{\operatorname{T}}\mathcal{A} \leq m$.
	\end{enumerate}
	Then, for any minimal component $(K_{0}, \varphi)$ of $(\mathcal{X},\varphi)$, we have that
	\begin{equation}
		\dim_{\operatorname{T}}(K_{0},\varphi) = \operatorname{rank}\mathfrak{M}_{\mathbb{Z}}(K_{0},\varphi) \leq m-1.
	\end{equation}
\end{theorem}
\begin{proof}
	Let $\pi_{\mathcal{X}}$ denote the restriction of $\pi$ to $\mathcal{X}$, which is a flow on $\mathcal{X}$ according to Theorem \ref{TH: CartwrightTheoremFlows}. Next, we apply Theorem \ref{TH: BendixsonCocycleContTime} with $j_{s} \coloneq \pi^{s}_{\mathcal{X}} \colon \Gamma \to \mathcal{B}_{s}$ and $\mathcal{A}$, $\Gamma$, and $\{\mathcal{B}_{s}\}_{s \leq 0}$ as in the statement. According Remark \ref{REM: BendixsonForContTimeCocyclesRemark}, all conditions of the referenced theorem are satisfied. Thus, $\check{H}^{m}(\Gamma;\mathbb{Z}) = 0$. As in Theorem \ref{TH: CartwrightTheoremSkewProductFlows}, $\Gamma$ can be endowed with the structure of a compact abelian group such that $K_{0}$ becomes its identity component. Thus, similarly to Theorem \ref{TH: CartwrightRefinementFlows}, we must have that $\dim_{\operatorname{T}}K_{0} = \dim_{\operatorname{T}} \Gamma \leq m-1$.
\end{proof}

In particular, Theorem \ref{TH: CartwrightRefinementSkewProduct} with $m=1$ guarantees that the minimal components of the vertical flow are stationary. However, it says nothing about the zero-dimensional topology of the fibers determined by torsion elements in the quotient $\mathfrak{M}_{\mathbb{Z}}(\pi)/\mathfrak{M}_{\mathbb{Z}}(\vartheta)$, which can be quite complex. This can be illustrated by the standard example of a solenoid, which arises as the inverse limit of a sequence of extensions with resonances corresponding to period doubling, as in the following example.
\begin{example}
Let $\mathbb{S}_{n} \coloneq \mathbb{R}/2^{n}\mathbb{Z}$ be the circle of length $2^{n}$. Define $\pi^{t}_{n}(x) \coloneq x + t \pmod {2^{n}}$ for $t \in \mathbb{R}$ and $x \in \mathbb{S}_{n}$. Then there is an extension map $\wp_{n} \colon (\mathbb{S}_{n+1},\pi_{n+1}) \to (\mathbb{S}_{n},\pi_{n})$ given by $\wp_{n}(x) = x \pmod {2^{n}}$, which is a two-to-one covering. Clearly, $\mathfrak{M}_{\mathbb{Z}}(\pi_{n}) = 2^{-n}\mathbb{Z}$, and $\mathfrak{M}_{\mathbb{Z}}(\pi_{n})/\mathfrak{M}_{\mathbb{Z}}(\pi_{0}) = 2^{-n}\mathbb{Z}/\mathbb{Z} = \operatorname{Tor}(\pi_{n},\pi_{0})$. The solenoid $(\mathbb{S},\pi)$ is obtained as the inverse limit of such systems:
\begin{equation}
	(\mathbb{S},\pi) \coloneq \varprojlim (\mathbb{S}_{n},\pi_{n}).
\end{equation}
By general arguments, see Theorem \ref{TH: AlmostPeriodicInverseLimit} below, it can be shown that $(\mathbb{S},\pi)$ is a minimal almost periodic flow with $\mathfrak{M}_{\mathbb{Z}}(\pi) = \mathbb{Z}[1/2] = \bigcup_{n \geq 0} \mathfrak{M}_{\mathbb{Z}}(\pi_{n})$. Considering $(\mathbb{S},\pi)$ as an extension of $(\mathbb{S}_{0},\pi_{0})$ via the mapping $\mathbb{S} \ni x=(x_{0},x_{1},\ldots) \overset{\wp}{\to} x_{0} \in \mathbb{S}_{0}$, we have that $\mathfrak{M}_{\mathbb{Z}}(\pi)/\mathfrak{M}_{\mathbb{Z}}(\pi_{0}) = \operatorname{Tor}(\pi,\pi_{0}) = \bigcup_{n \geq 0} 2^{-n}\mathbb{Z}/\mathbb{Z}$ and
\begin{equation}
	\operatorname{Hom}(\operatorname{Tor}(\pi,\pi_{0});\mathbb{T}^{1}) = \varprojlim \operatorname{Hom}(\operatorname{Tor}(\pi_{n},\pi_{0});\mathbb{T}^{1}) = \varprojlim \mathbb{Z}/2^{n}\mathbb{Z} = \mathbb{Z}_{2},
\end{equation}
where $\mathbb{Z}_{2}$ is the group of $2$-adic integers. Since the minimal components of $\pi$ are just singletons, the fiber over $0 \in \mathbb{S}_{0}$ is isomorphic to $\mathbb{Z}_{2}$, which is topologically a Cantor set. Note that $\mathbb{S}$ cannot be homeomorphic to $\mathbb{Z}_{2} \times \mathbb{S}_{0}$ since it is connected, whereas the product has uncountably many connected components. In other words, $\wp \colon \mathbb{S} \to \mathbb{S}_{0}$ is only a locally trivial bundle over $\mathbb{S}_{0}$ with fiber $\mathbb{Z}_{2}$. In terms of $\mathbb{Z}_{2}$, the monodromy mapping $\pi^{1} \colon \mathbb{Z}_{2} \to \mathbb{Z}_{2}$ is the dyadic odometer: $a \mapsto a + 1$. \qed
\end{example}

For general inverse limits of minimal almost periodic flows, one can state the following. 
\begin{theorem}
	\label{TH: AlmostPeriodicInverseLimit}
	Let $(\mathcal{X},\pi)$ be the inverse limit of a sequence of extensions $\wp_{n} \colon (\mathcal{X}_{n+1},\pi_{n+1}) \to (\mathcal{X}_{n},\pi_{n})$, where $n=0,1,\ldots$, of minimal almost periodic flows. Then $(\mathcal{X},\pi)$ is a minimal almost periodic flow with $\mathfrak{M}_{\mathbb{Z}}(\pi) = \bigcup_{n \geq 0} \mathfrak{M}_{\mathbb{Z}}(\pi_{n})$. Moreover,
	\begin{enumerate}
		\item[1)] $\mathfrak{M}_{\mathbb{Z}}(\pi)/\mathfrak{M}_{\mathbb{Z}}(\pi_{0}) = \bigcup_{n \geq 0} \mathfrak{M}_{\mathbb{Z}}(\pi_{n})/\mathfrak{M}_{\mathbb{Z}}(\pi_{0})$ and $\operatorname{Tor}(\pi, \pi_{0}) = \bigcup_{n \geq 0} \operatorname{Tor}(\pi_{n},\pi_{0})$;
		\item[2)] $\operatorname{Hom}(\operatorname{Tor}(\pi,\pi_{0});\mathbb{T}) = \varprojlim \operatorname{Hom}(\operatorname{Tor}(\pi_{n},\pi_{0});\mathbb{T})$;
		\item[3)] One can construct vertical flows $\varphi_{n}$ of the extensions $\wp_{n}$ such that the inverse limit $\varphi$ of $\wp_{n} \colon (\mathcal{X}_{n+1}, \varphi_{n+1}) \to (\mathcal{X}_{n},\varphi_{n})$ is a vertical flow of $\wp \colon (\mathcal{X},\pi) \to (\mathcal{X}_{0},\pi_{0})$.
	\end{enumerate}
\end{theorem}
\begin{proof}
	Fix $x=(x_{0},x_{1},\ldots) \in \mathcal{X}$. By definition, $x_{n} \in \mathcal{X}_{n}$ and $\wp_{n}(x_{n+1}) = x_{n}$ for any $n$. Let $G_{n}$ be a $\pi_{n}$-consistent group structure on $\mathcal{X}_{n}$ with neutral element $x_{n}$. Then $\wp_{n} \colon G_{n+1} \to G_{n}$ becomes a surjective homomorphism. Let $G$ be the inverse limit of $G_{n}$. By general properties, $G$ is a metrizable compact connected abelian group, for which $x$ is the neutral element. It is easy to verify that $G$ is $\pi$-consistent and that the orbit of $x$ is dense in $\mathcal{X}$. Consequently, $\pi$ is a minimal almost periodic flow.
	
	Since $\wp_{n}$ is surjective, it induces an injective homomorphism on the character groups: $\wp^{*}_{n} \colon \hat{G}_{n} \to \hat{G}_{n+1}$. By duality, $\hat{G} = \varinjlim \hat{G}_{n}$. Since frequency modules can be expressed in terms of character groups using homomorphisms $i_{\pi_{n}}(t) = \pi^{t}_{n}(x_{n})$, which, given our choice of points, agree with $\wp_{n}$, items 1) and 2) become obvious.
	
	For item 3), we take $\varphi^{t}_{n}(g) = g + \pi^{t}_{n}(x_{n})$ for $g \in G_{n}$ and $\varphi^{t}(g) = g + \pi^{t}(x)$ for $g \in G$. Then it is straightforward to verify that these are the required vertical flows.
\end{proof}

In concluding this section, we note that standard approaches to establishing the existence of almost periodic solutions a priori assume positive uniform Lyapunov stability, see \cite{LevitanZhikov1982, SackerSell2017, Cheban2020, ShenYi1998}. This, in particular, guarantees distality and is, indeed, a necessary condition, since almost periodic flows are precisely those flows that possess two-sided uniform Lyapunov stability. It turns out that this a priori assumption is sufficient under certain constraints on dimension or order preservation. However, with the exception of \cite[Theorem 5, Chapter 7, Section 8]{LevitanZhikov1982}, all these results correspond to extensions with trivial torsion-free part of the quotient $\mathfrak{M}_{\mathbb{Z}}(\pi)/\mathfrak{M}_{\mathbb{Z}}(\vartheta)$. From this perspective, Cartwright’s results are interesting in that they can be applied in a more general context, even if suitable existence theorems have not been obtained.

A step beyond this is to consider minimal almost automorphic flows $(\mathcal{X},\pi)$ that extend a minimal almost periodic flow $(\mathcal{Q},\vartheta)$. As shown in \cite{ShenYi1998}, such a situation arises in a number of cases if the requirement of uniform Lyapunov stability is dropped. Moreover, almost automorphy is closely related to so-called strange nonchaotic attractors, which have been observed experimentally in systems subject to quasiperiodic external forcing (see \cite{AnikushinAADyn2021} and references therein).

Any minimal almost automorphic flow $(\mathcal{X},\pi)$ is an almost 1-1 extension\footnote{That is, the set of $y \in \mathcal{Y}$ whose preimage is a singleton is residual.} of a minimal almost periodic flow $(\mathcal{Y},\psi)$, which is unique up to a topological conjugacy and is called the \textit{maximal almost periodic factor} of $(\mathcal{X},\pi)$. Specifically, $\mathcal{Y}$ is obtained by taking the quotient of $\mathcal{X}$ with respect to the proximal relation, see \cite[Chapter I, Corollary 2.15]{ShenYi1998}. If $\wp_{\mathcal{Y}} \colon \mathcal{X} \to \mathcal{Y}$ is the projection onto the quotient, then $\wp$ factors through $\wp_{\mathcal{Y}}$, yielding
\begin{equation}
	\label{EQ: AlmostAutomorphicTriple}
	(\mathcal{X},\pi) \overset{\wp_{\mathcal{Y}}}{\to} (\mathcal{Y},\psi) \overset{\wp_{\mathcal{Q}}}{\to} (\mathcal{Q},\vartheta)
\end{equation}
with $\wp = \wp_{\mathcal{Y}} \circ \wp_{Q}$. The extension $\wp_{\mathcal{Q}}$ is fully described by Theorem \ref{TH: CartwrightTheoremSkewProductFlows}. It would be interesting to obtain bounds for $\mathfrak{M}_{\mathbb{Q}}(\psi)/\mathfrak{M}_{\mathbb{Q}}(\vartheta)$ using certain constraints for $(\mathcal{X},\pi)$ in the spirit of Theorem \ref{TH: CartwrightRefinementSkewProduct}. However, in the general case the situation turns out to be complicated, since the projection $\wp_{\mathcal{Y}}$ may increase the dimension of the fiber $\mathcal{Y}_{q} = \wp_{\mathcal{Y}}(\mathcal{X}_{q})$. The following example illustrates this. 
\begin{example}
\label{EX: DenjoyAlmostAutomorphic}
Let $f$ be a Denjoy homeomorphism of the circle $\mathbb{S}^{1}=\mathbb{R}/\mathbb{Z}$ with irrational rotation number $\omega$, and let $h \colon \mathbb{S}^{1} \to \mathbb{S}^{1}$ provide the extension of the irrational rotation $(\mathbb{S}^{1},R_{\omega})$ to $(\mathbb{S}^{1},f)$. Recall that, according to the construction, we choose an orbit $x_{k} = R^{k}_{\omega}(x_{0})$ of $R_{\omega}$ and construct a new circle for $f$ by blowing intervals\footnote{The intervals are chosen in such a way that the series of lengths $|\mathcal{I}_{k}|$ is convergent.} $\mathcal{I}_{k}=[a_{k},b_{k}]$ around each point $x_{k}$, where $k \in \mathbb{Z}$. Then $f$ linearly maps $\mathcal{I}_{k}$ onto $\mathcal{I}_{k+1}$ and acts as $R_{\omega}$ on the complement, whereas the extension $h$ shrinks each interval $\mathcal{I}_{k}$ back to $x_{k}$ and is identical on their complement. Then there is a unique minimal Cantor-like set $\mathcal{C} = \mathbb{S}^{1} \setminus \bigcup_{k \in \mathbb{Z}}(a_{k},b_{k})$ of $f$. Define $(\mathcal{X},\pi)$ as the mapping torus 
\begin{equation}
	\mathcal{X}=\mathcal{C} \times [0,1]/((x,1) \sim (f(x),0))
\end{equation}
over $f$ restricted to $\mathcal{C}$, and let $(\mathcal{Q},\vartheta)$ be the linear flow on $\mathbb{S}^{1}$: $\vartheta^{t}(\theta) = \theta + t \pmod 1$. Then $\wp \colon \mathcal{X} \to \mathbb{S}^{1}$ given by $(x,\theta) \mapsto \theta$ is an extension of $(\mathcal{Q}, \vartheta)$ to $(\mathcal{X},\pi)$. Finally, let $(\mathcal{Y},\psi)$ be the linear flow $(\theta_{1},\theta_{2}) \mapsto (\theta_{1} + \omega t, \theta_{2} + t)$ on the $2$-torus $\mathbb{T}^{2} = \mathbb{R}^{2}/\mathbb{Z}^{2}$, which is the mapping torus over $(\mathbb{S}^{1},R_{\omega})$. With $\wp_{\mathcal{Y}}(x,\theta) = (h(x),\theta)$ and $\wp_{\mathcal{Q}}(\theta_{1},\theta_{2}) = \theta_{2}$, we are in the situation of \eqref{EQ: AlmostAutomorphicTriple}. When restricting $h$ to $\mathcal{C}$, the injectivity is violated only at the endpoints $\{a_{k}, b_{k}\}$, which are glued together into $x_{k}$. It is not hard to see that the corresponding orbits in $\mathcal{X}$ are mapped by $\wp_{\mathcal{Y}}$ onto the entire orbit $\mathcal{Y}^{\operatorname{sing}}$ of $(x_{0},0)$. For $y$ in the residual set $\mathcal{Y}^{\operatorname{reg}} \coloneq \mathcal{Y} \setminus \mathcal{Y}^{\operatorname{sing}}$, the preimage $\wp^{-1}_{\mathcal{Y}}(y)$ is a singleton. By \cite[Theorem 2.14, Chapter I]{ShenYi1998}, $(\mathcal{X},\pi)$ is a minimal almost automorphic flow and $(\mathcal{Y},\psi)$ is its maximal almost periodic factor. Clearly, for any $q \in \mathcal{Q}$, the fiber $\mathcal{X}_{q} = \wp^{-1}(q)$ is a Cantor set, while $\mathcal{Y}_{q} = \wp^{-1}_{\mathcal{Y}}(q)$ is a circle. \qed
\end{example}

In \cite{ShenYi1998}, the constraints regarding dimension and ordering rule out the phenomena considered in the above example: $(\mathcal{Y},\psi)$ turns out to be an $N$-to-1 extension of $(\mathcal{Q},\vartheta)$ for some $N \geq 1$. We find it interesting to attempt to obtain estimates for frequency modules in the general case, even in the absence of existence theorems. To illustrate the difficulties, let $\mathcal{Y}^{\operatorname{reg}}$ be the set of $y \in \mathcal{Y}$ such that $\wp^{-1}_{\mathcal{Y}}(y)$ is a singleton. Then $\wp^{-1}_{\mathcal{Y}}$ is continuous on $\mathcal{Y}^{\operatorname{reg}}$, see \cite[Corollaries 2.15 and 2.17, Chapter I]{ShenYi1998}. For any $q \in \mathcal{Q}$, the fiber $\mathcal{Y}_{q}$ splits into the regular part $\mathcal{Y}^{\operatorname{reg}}_{q} \coloneq \mathcal{Y}_{q} \cap \mathcal{Y}^{\operatorname{reg}}$ and the singular part $\mathcal{Y}^{\operatorname{sing}}_{q} \coloneq \mathcal{Y}_{q} \setminus \mathcal{Y}^{\operatorname{reg}}_{q}$. Since $\mathcal{Y}^{\operatorname{reg}}_{q}$ is homeomorphic to a subspace of $\mathcal{X}_{q}$ via $\wp_{\mathcal{Y}}$, the addition theorem (see \cite[Theorem 1.5.10]{Engelking1978}) ensures that
\begin{equation}
	\dim_{\operatorname{T}}\mathcal{Y}_{q} \leq \dim_{\operatorname{T}} \mathcal{X}_{q} + \dim_{\operatorname{T}} \mathcal{Y}^{\operatorname{sing}}_{q} + 1.
\end{equation}
In Example \ref{EX: DenjoyAlmostAutomorphic}, we have $\dim_{\operatorname{T}} \mathcal{X}_{q} = 0$ and $\dim_{\operatorname{T}} \mathcal{Y}^{\operatorname{sing}}_{q} = 0$, since $\mathcal{Y}^{\operatorname{sing}}_{q}$ is countable. Thus, the estimate is sharp in general. It is also clear that the argument similar to Theorem \ref{TH: CartwrightRefinementSkewProduct} can be applied in the almost automorphic case only when $\wp^{*}_{\mathcal{Y}} \colon \check{H}^{m}(\mathcal{Y}_{q};\mathbb{Z}) \to \check{H}^{m}(\mathcal{X}_{q};\mathbb{Z})$ is injective, that is violated in the example.




\section*{Acknowledgments}
The author is grateful to Claude Sonnet 5.5 for fruitful discussions.

\begin{funding}
The reported study was funded by the Russian Science Foundation (Project 25-11-00147).
\end{funding}

\section*{Data availability}
Data sharing not applicable to this article as no datasets were generated or analyzed during the current study.

\section*{Conflict of interest}
The authors declare that they have no conflict of interest.



\begin{thebibliography}{99}







\bibitem{AnikushinRomanovQBC2026}
Anikushin~M.M., Romanov~A.O. Quadratic Bendixson criterion for delay systems. \textit{arXiv preprint} (2026)

\bibitem{AnikushinRomanov2026Subroutine}
Anikushin~M.M., Romanov~A.O. Quadratic Bendixson criterion for delay systems: a numerical implementation. \textit{arXiv preprint} (2026)

\bibitem{Anikushin2023Comp}
Anikushin~M.M. Spectral comparison of compound cocycles generated by delay equations in Hilbert spaces. \textit{arXiv preprint}, arXiv:2302.02537 (2026)

\bibitem{AnikushinRomanov2023FreqConds}
Anikushin~M.M., Romanov~A.O. Frequency conditions for the global stability of nonlinear delay equations with several equilibria. \textit{J. Math. Anal. Appl.}, \textbf{563}(2), 130775 (2026)

\bibitem{Anikushin2023LyapExp}
Anikushin~M.M. Variational description of uniform Lyapunov exponents via adapted metrics on exterior products. \textit{Nonlinear Differ. Equ. Appl.}, \textbf{33}(1), 26 (2026)

\bibitem{AnikushinRomanov2024EffEst}
Anikushin~M.M., Romanov~A.O. On the recent progress in effective dimension estimates for delay equations. \textit{Differ. Uravn. Protsessy Upravl.}, 1, 22--46 (2024)

\bibitem{Anikushin2020FreqDelay}
Anikushin~M.M. Frequency theorem and inertial manifolds for neutral delay equations. \textit{J. Evol. Equ.}, \textbf{23}, 66 (2023)



\bibitem{Anikushin2020FreqParab}
Anikushin M.M. Frequency theorem for parabolic equations and its relation to inertial manifolds theory. \textit{J. Math. Anal. Appl.}, \textbf{505}(1), 125454 (2021)

\bibitem{AnikushinAADyn2021}
Anikushin~M.M. Almost automorphic dynamics in almost periodic cocycles with one-dimensional inertial manifold. \textit{Differ. Uravn. Protsessy Upravl.}, 2, (2021) (in Russian)

\bibitem{Anikushin2019Liouv}
Anikushin~M.M. On the Liouville phenomenon in estimates of fractal dimensions of forced quasi-periodic oscillations. \textit{Vestn. St. Petersbg. Univ., Math.}, \textbf{52}(3), 234--243 (2019)

\bibitem{AnikushinReitRom2019}
Anikushin~M.M., Reitmann V., Romanov A.O. Analytical and numerical estimates of the fractal dimension of forced quasiperiodic oscillations in control systems. \textit{Differ. Uravn. Protsessy Upravl.}, 2 (2019), (in Russian)

\bibitem{BargeGiraldoSanjurjo2020}
Barge~H., Giraldo~A., Sanjurjo~J.M. Bifurcations, robustness and shape of attractors of discrete dynamical systems. \textit{J. Fixed Point Theory Appl.}, \textbf{22}(2), 29 (2020)

\bibitem{Bendixson1901}
Bendixson~I.O. Sur les courbes d\'{e}finies par des \'{e}quations diff\'{e}rentielles. \textit{Acta Math.}, \textbf{24}, 1--88 (1901)


\bibitem{Cartwright1969}
Cartwright~M.L. Almost periodic differential equations and almost periodic flows. \textit{J. Differ. Equ.}, \textbf{5}(1), 167--181 (1969)

\bibitem{Cartwright1967}
Cartwright~M.L. Almost periodic flows and solutions of differential equations. \textit{Proc. Lond. Math. Soc.}, \textbf{3}(4), 768--768 (1967)

\bibitem{Cheban2020}
Cheban~D.N. \textit{Nonautonomous Dynamics. Nonlinear Oscillations and Global Attractors}, Springer (2020)

\bibitem{EilenbergSteenrod2015}
Eilenberg~S., Steenrod~N. \textit{Foundations of Algebraic Topology}. Princeton University Press (1952)

\bibitem{Engelking1978}
Engelking~R. \textit{Dimension Theory}. Amsterdam: North-Holland Publishing Company (1978)

\bibitem{GutmanJinTsukamoto2019}
Gutman~Y., Jin~L., Tsukamoto~M. A Lipschitz refinement of the Bebutov--Kakutani dynamical embedding theorem. \textit{J. Dyn. Differ. Equ.}, \textbf{31}(1), 81--91 (2019)

\bibitem{HeZhouLv2002}
He~L., Zhou~L., Lv~J. The dynamical properties of a continuous semiflow and its inverse limit. \textit{Chaos Solit.}, \textbf{14}(8), 1169--1177 (2002)

\bibitem{HofmannMostert1973}
Hofmann~K.H., Mostert~P.S. \textit{Cohomology Theories for Compact Abelian Groups}. Springer (1973)

\bibitem{KapitanskiRodnianski2000}
Kapitanski~L., Rodnianski~I. Shape and Morse theory of attractors. \textit{Commun. Pure Appl. Math.}, \textbf{53}(2), 218--242 (2000)

\bibitem{Kuzetal2020Lorenz}
Kuznetsov~N.V., Mokaev~T.N., Kuznetsova~O.A., Kudryashova~E.V. The Lorenz system: hidden boundary of practical stability and the Lyapunov dimension. \textit{Nonlinear Dyn.}, \textbf{102}, 713--732 (2020)

\bibitem{KuzReit2020}
Kuznetsov~N.V., Reitmann~V. \textit{Attractor Dimension Estimates for Dynamical Systems: Theory and Computation}. Switzerland: Springer International Publishing AG (2020)


\bibitem{LeonovBoi1992}
Leonov~G.A., Boichenko~V.A. Lyapunov's direct method in the estimation of the Hausdorff dimension of attractors. \textit{Acta Appl. Math.}, \textbf{26}(1), 1--60 (1992)

\bibitem{LevitanZhikov1982}
Levitan~B.M., Zhikov~V.V. \textit{Almost Periodic Functions and Differential Equations}. CUP Archive (1982)

\bibitem{LiMuldowney1996SIAMGlobStab}
Li~M.Y., Muldowney~J.S. A geometric approach to global-stability problems. \textit{SIAM J. Math. Anal.}, \textbf{27}(4), 1070--1083 (1996)

\bibitem{LiMuldowney1995LowBounds}
Li~M.Y., Muldowney~J.S. Lower bounds for the Hausdorff dimension of attractors. \textit{J. Dyn. Differ. Equ.}, \textbf{7}(3), 457--469 (1995)

\bibitem{MalletParet1976}
Mallet-Paret~J. Negatively invariant sets of compact maps and an extension of a theorem of Cartwright. \textit{J. Differ. Equ.}, \textbf{22}(2), 331--348 (1976)

\bibitem{Morita1975}
Morita~K. \v{C}ech cohomology and covering dimension for topological spaces. \textit{Fundam. Math.}, 87, 31--52 (1975)

\bibitem{NakajimaSeifert1983}
Nakajima~F., Seifert~G. The number of periodic solutions of 2-dimensional periodic systems. \textit{J. Differ. Equ.}, \textbf{49}(3), 430--440 (1983)

\bibitem{NoackReitmann1996}
Noack~A., Reitmann~V. Hausdorff dimension estimates for invariant sets of time-dependent vector fields. \textit{Z. Anal. Anwend.}, \textbf{15}(2), 457--473 (1996)


\bibitem{Pontryagin1973}
Pontryagin~L.S. \textit{Topological Groups}. Moscow (1973) (in Russian)

\bibitem{SackerSell2017}
Sacker~R.J., Sell~G.R. Almost periodicity, Ricker map, Beverton--Holt map and others, a general method. \textit{J. Differ. Equ. Appl.}, \textbf{23}(7), 1286--1297 (2017)

\bibitem{ShenYi1998}
Shen~W., Yi~Y. Almost Automorphic and Almost Periodic Dynamics in Skew-Product Semiflows. AMS (1998)

\bibitem{Smith1986HD}
Smith~R.A. Some applications of Hausdorff dimension inequalities for ordinary differential equations. \textit{P. Roy. Soc. Edinb. A}, \textbf{104}(3-4), 235--259 (1986)

\bibitem{Spanier1966}
Spanier~E.H. \textit{Algebraic Topology}. Springer (1966)

\bibitem{Thieullen1992}
Thieullen~P. Entropy and the Hausdorff dimension for infinite-dimensional dynamical systems. \textit{J. Dynam. Differ. Equ.}, \textbf{4}, 127--159 (1992)


\bibitem{ZelikAttractors2022}
Zelik~S. Attractors. Then and now. \textit{Russ. Math. Surv.} \textbf{78}(4), 635--777 (2023) and \textit{Uspekhi Mat. Nauk}, \textbf{78}(4), 53--198 (2023)

\bibitem{Zinchenko1988}
Zinchenko~I.L. The group of characters on the closure of the almost periodic trajectory of an autonomous system of differential equations. \textit{Diff. Urav.}, \textbf{24}(6), 1043--1045 (1988) (in Russian)
\end{thebibliography}
\end{document}